\documentclass[12pt,a4paper]{article}
\usepackage{amsmath}
\usepackage[font={footnotesize,it}]{caption}
\usepackage[top=1.2in, bottom=1.2in, left=1.1in, right=1.1in]{geometry}

\DeclareCaptionStyle{italic}[justification=centering]{labelfont={bf},textfont={it},labelsep=colon}
\usepackage{graphicx,psfrag,epsf}
\usepackage{enumerate}
\usepackage{natbib}
\usepackage{url, setspace}
\usepackage{booktabs, subfig, bm, paralist,mathpazo,tikz,todonotes,longtable,microtype,dsfont,rotating} 
\usepackage[pdftex,colorlinks=true]{hyperref}
\definecolor{darkblue}{rgb}{0,0,.6}
\hypersetup{citecolor=darkblue,linkcolor=darkblue,urlcolor=darkblue}

\newcommand{\blind}{0}

\newcommand{\X}{\mathcal{X}}

\graphicspath{{plots/}}
\DeclareMathOperator*{\argmin}{\arg\!\min}
\newsavebox\CBox
\def\textBF#1{\sbox\CBox{#1}\resizebox{\wd\CBox}{\ht\CBox}{\textbf{#1}}}

\date{\today}
\begin{document}

\def\spacingset#1{\renewcommand{\baselinestretch}
{#1}\small\normalsize} \spacingset{1}

\if0\blind
{
  \title{\bf A comparison of Hurst exponent estimators in long-range dependent curve time series}
  \author{Han Lin Shang\thanks{Postal address: Research School of Finance, Actuarial Studies and Statistics, Level 4, Building 26C, Australian National University, Kingsley Street, Acton, Canberra, ACT 2601, Australia; Telephone: +61(2) 612 50535; Fax: +61(2) 612 50087; Email: hanlin.shang@anu.edu.au; ORCID ID: \url{https://orcid.org/0000-0003-1769-6430}.}
  \hspace{.2cm}\\
    Research School of Finance, Actuarial Studies and Statistics \\
    Australian National University
}
  \maketitle
} \fi

\if1\blind
{
   \title{\bf A comparison of Hurst exponent estimators in long-range dependent curve time series}
   \author{}
   \maketitle
} \fi

\bigskip
\begin{abstract}

The Hurst exponent is the simplest numerical summary of self-similar long-range dependent stochastic processes. We consider the estimation of Hurst exponent in long-range dependent curve time series. Our estimation method begins by constructing an estimate of the long-run covariance function, which we use, via dynamic functional principal component analysis, in estimating the orthonormal functions spanning the dominant sub-space of functional time series. Within the context of functional autoregressive fractionally integrated moving average models, we compare finite-sample bias, variance and mean square error among some time- and frequency-domain Hurst exponent estimators and make our recommendations.

\vspace{.1in}
\end{abstract}

\newpage

\def\spacingset#1{\renewcommand{\baselinestretch}{#1}\small\normalsize} \spacingset{1}
\spacingset{1.45}

\section{Introduction}\label{sec:1}

In univariate time series analysis, long memory was brought to prominence by \cite{Hurst51} and \cite{Mandelbrot63}, and it has subsequently received extensive attention in the literature \citep[see, e.g.,][]{Beran94, EM02, DOT03, Robinson03, Palma07}. Of importance in analyzing and modeling long-memory univariate time series is estimating the strength of the long-memory dependence. There are two measures commonly used: The parameter $H$, known as the Hurst exponent or self-similarity parameter \citep{MV68} and the fractional integration parameter, $d$, arises from the generalization of autoregressive fractionally integrated moving average (ARFIMA$(p, d, q)$) models from integer to non-integer values of the integration parameter $d$. The two parameters are closely related through the simple formula $H=d+\frac{1}{2}$.

In univariate time series analysis, a number of Hurst exponent estimators have been developed, and theoretical results on the asymptotic properties of various estimators have been obtained. Because the finite-sample properties of these estimators can be quite different from their asymptotic properties, several authors considered an empirical comparison of estimators of $H$ and $d$. Nine estimators were discussed in some detail by \cite{TTW95} who performed an empirical investigation of these estimators for a single series length of 10,000 data points, five values of both $H$ and $d$, and 50 replications. \cite{TT97} showed in a simulation study that the differenced variance estimator was unbiased for five values of $H$ (0.5, 0.6, 0.7, 0.8 and 0.9) for series with 10,000 observations whereas the aggregated variance estimator was downwards biased. \cite{Jensen99} presented a comparison of two estimators based on wavelets, and a Geweke-Porter-Hudak (GPH) estimator for four series lengths ($2^7, 2^8, 2^9, 2^{10}$ observations), five values of $d$ and 1,000 replications. \cite{JLM+07} performed a comparison of six estimators on simulated fractional Gaussian noise with $2^{15}$ observations, five values of $H$ and 100 replications. 

Long-memory functional time series analysis was recently studied by \cite{LRS18}, who proposed an R/S estimation method for determining long-memory parameter in a functional ARFIMA model, where observations are temporally dependent continuous functions, for example, age-specific fertility rate improvement observed over the years \citep[e.g.,][]{HU07, CM09}. The functional ARFIMA model can be viewed as a generalization of many parametric models. For example, \cite{Bosq00} and \cite{BB07} provided the functional autoregressive of order 1 (FAR(1)) and derived one-step-ahead forecasts that are based on a regularized form of the Yule-Walker equations. Later, FAR(1) was extended to FAR($p$), where the order $p$ can be determined via a sequential hypothesis testing procedure \citep{KR13}. \cite{ANH15} proposed a forecasting method based on vector autoregressive (VAR) forecasts of principal component scores. The method of \cite{ANH15} can also be viewed as an extension of \cite{HS09}, where principal component scores are forecast via a univariate time series forecasting method. \cite{KK17} considered the functional moving average (FMA) process and introduced an innovation algorithm to obtain the best linear predictor. \cite{KKW17} extended the VAR model to vector autoregressive moving average model for modeling and forecasting principal component scores, which can be viewed as a simpler estimation approach of the functional autoregressive moving average. \cite{AK17} showed the equivalent relationship between FMA and vector moving average. 

A central issue in functional time series analysis is to model the temporal dependence of the functional observations accurately. Following the early work of \cite{LRS18}, we compare the finite-sample estimation accuracy of several Hurst exponent estimators in functional ARFIMA models. Our method constructs an estimate of the long-run covariance function, which we use, via dynamic functional principal component analysis, in estimating the orthonormal functions spanning the dominant sub-space of the curves. Based on the first set of principal component scores, we apply several univariate time series Hurst exponent estimators, and compare their estimation accuracy in terms of bias, variance and mean square error. Our goal is to provide some practical guidance on the method that provides the best estimation accuracy of the Hurst exponent.

The remainder of the paper is outlined as follows. In Section~\ref{sec:2}, we present two methods for estimating long-run covariance function, from which the dominant set of principal component scores can be obtained. In Section~\ref{sec:3}, we revisit some long-memory univariate time series estimators for estimating the Hurst exponent. In Section~\ref{sec:4}, we compare the estimation accuracy of various estimators and make our recommendation. 

\section{Dynamic functional principal component analysis}\label{sec:2}

\subsection{Estimation of the long-run covariance function}

A time series of functions can be denoted as $\{\X_t, t\in \mathbb{Z}\}$, where $\mathbb{Z}=\{t: t=0,\pm 1, \dots\}$ and each $\X_t$ is a random function of a stochastic process $\X(\omega)$ where $\omega\in \mathcal{I}\subset R$ represents a continuum bounded within a finite interval of the real line. Further, let $\{\X_t(\omega)\}_{t\in\mathbb{Z}}$ be a stationary and ergodic functional time series. For a stationary functional time series, the long-run covariance operator is defined as
\begin{align*}
C(\omega, \nu) &= \sum_{\ell=-\infty}^{\infty}\gamma_{\ell}(\omega, \nu) \\
\gamma_{\ell}(\omega, \nu) &= \text{cov}[\X_0(\omega), \X_{\ell}(\nu)],
\end{align*}
and is a well-defined element of $\mathcal{L}^2(\mathcal{I})^2$ for a compact support interval $\mathcal{I}$, under mild weak dependence and moment conditions. By assuming $\X$ is a continuous and square-integrable function, the function $\mathcal{C}$ induces the kernel operator $\mathcal{L}^2(\mathcal{I})\mapsto \mathcal{L}^2(\mathcal{I})$. Through right integration, $\mathcal{C}$ defines a Hilbert-Schmidt integral operator on $\mathcal{L}^2(\mathcal{I})$ given by
\begin{equation*}
\mathcal{C}(\X)(\omega) = \int_{\mathcal{I}}C(\omega, \nu)\X(\nu)d\nu,
\end{equation*}
whose eigenvalues and eigenfunctions are related to the dynamic functional principal components defined in \cite{HKH15}. 

In practice, we need estimate $C(\omega, \nu)$ from a finite sample $\bm{\X}(\omega) = \left[\X_1(\omega), \dots, \X_n(\omega)\right]$. Given its definition as a bi-infinite sum, a natural estimator of $C$ is
\begin{equation}
\widehat{C}_n(\omega, \nu) = \frac{1}{n^{3-2\alpha}}\sum_{|\ell|=0}^{|\ell|\leq n}\left(n - |\ell|\right)\widehat{\gamma}_{\ell}(\omega, \nu),  \label{eq:long-cov}
\end{equation}
where $\alpha=\frac{3}{2}-H$ is the so-called memory parameter, $\ell$ denotes a lag variable, and 
\begin{equation*}
\widehat{\gamma}_{\ell}(\omega, \nu) = \left\{ \begin{array}{ll}
         \frac{1}{n}\sum_{j=1}^{n-\ell}\left[\X_j(\omega) - \overline{\X}(\omega)\right]\left[\X_{j+\ell}(\nu) - \overline{\X}(\nu)\right] & \mbox{if $\ell \geq 0$};\\
        \frac{1}{n}\sum^n_{j=1-\ell}\left[\X_j(\omega) - \overline{\X}(\omega)\right]\left[\X_{j+\ell}(\nu) - \overline{\X}(\nu)\right] & \mbox{if $\ell < 0$}.\end{array} \right. 
\end{equation*}        
is an estimator of $\gamma_{\ell}(\omega, \nu)$. In the case of stationary short-memory functional time series, it is known that $\alpha=1$. From~\eqref{eq:long-cov}, the estimated long-run covariance is obtained by summing all autocovariance functions with linearly decreasing weights. Let $W$ denote the number of grid points in a curve. In \cite{LRS18}, they consider $\ell=\min(n, W)$. For instance, when $W\geq n$, all finite-order lags are utilized. To estimate the value of $\alpha$, \cite{LRS18} applied the rescaled range (R/S) estimator of \cite{Hurst51} to the first set of dynamic principal component scores obtained from eigendecomposition of 
\begin{equation*}
\widehat{\widehat{C}}_n(\omega, \nu) = \sum_{|\ell|=0}^{|\ell|\leq n}\left(n - |\ell|\right)\widehat{\gamma}_{\ell}(\omega, \nu),
\end{equation*}
since $\frac{1}{n^{3-2\alpha}}$ in~\eqref{eq:long-cov} is a constant and it does not affect the estimation of the orthonormal functions spanning the dominant sub-space of functional time series.

\subsubsection{Kernel sandwich estimator}

Another long-run covariance estimator is the kernel sandwich estimator inspired by \cite{Andrews91} and \cite{AM92}. It is given by
\begin{equation}
\widehat{C}_{h,q}(\omega, \nu) = \sum_{\ell=-\infty}^{\infty}W_q\left(\frac{\ell}{h}\right)\widehat{\gamma}_{\ell}(\omega, \nu),\label{eq:1}
\end{equation}
where $h$ is called the bandwidth parameter and $W_q$ is a symmetric weight function with bounded support of order $q$. The kernel sandwich estimator in~\eqref{eq:1} was introduced in \cite{PT12}, \cite{HKR13}, \cite{RS17}, \citet[][Chapter 8.5]{KR17}, among others. As with any kernel estimator, the crucial part is on the estimation of bandwidth parameter $h$. It can be selected through a data-driven approach, such as the plug-in algorithm of \cite{RS17}. The plug-in bandwidth selection method can be summarized as:
\begin{enumerate}
\item[(1)] Compute pilot estimates of $C^{(p)}$, for $p=0$ and initial order of kernel function $q_1$:
\begin{equation*}
\widehat{C}_{h_1, q_1}^{(p)}(u, s) = \sum^{\infty}_{\ell = -\infty}W_{q_1}\left(\frac{\ell}{h_1}\right)|\ell|^p\widehat{\gamma}_{\ell}(u,s),
\end{equation*}
that utilize an initial bandwidth choice $h_1 = h_1(n)$, and weight function $W_{q_1}$ of order $q_1$.
\item[(2)] As established in \cite{BHR16}, estimate $c_0$ by
\begin{equation*}
\hspace{-.3in}{\widehat{c}_0(h_1, q_1, q) = (2q\|w\widehat{C}_{h_1, q_1}^{(q)}\|^2)^{1/(1+2q)}\left(\left(\|\widehat{C}_{h_1, q_1}^{(0)}\|^2+\left(\int^1_0\widehat{C}_{h_1, q_1}^{(0)}(u,u)du\right)^2\right)\int^{\infty}_{-\infty}W_{q_1}^2(x)dx\right)^{-1/(1+2q)},}
\end{equation*}
where $q$ denotes the final order of kernel function, $w$ is a constant depending on the final order of kernel function, and $\int^{\infty}_{-\infty}W_{q_1}^2(x)dx$ is a weight depending on the initial order of kernel function. A list of $w$ and $\int^{\infty}_{-\infty}W_{q_1}^2(x)dx$ values is presented in Table~\ref{tab:weights}.
\begin{table}[!htbp]
\centering
\tabcolsep 0.2in
\caption{A list of $w$ and $\int^{\infty}_{-\infty}W_{q_1}^2(x)dx$ values\label{tab:weights}}
\begin{tabular}{@{}lll@{}}
\toprule
Kernel function & $w$ & $\int^{\infty}_{-\infty}W_{q_1}^2(x)dx$ \\
\midrule
Bartlett & 1 & 2/3 \\
Parzen & 6 &  0.539285 \\
Tukey-Hanning & $\pi^2/4$ & 3/4 \\
Quadratic Spectral & $18\times \pi^2/125$ & 1 \\
Flat-top & & 4/3 \\
\bottomrule
\end{tabular}
\end{table}

For the initial kernel function, \cite{RS17} recommend to use flat-top kernel function, i.e., $q_1=\infty$. For the final kernel function, \cite{RS17} recommend to use Bartlett kernel function, i.e., $q=1$. Further, there exists $w$ satisfying $0<w=\lim_{x\rightarrow 0}|x|^{-q}(1-W_q(x))<\infty$.
\item[(3)] Use the bandwidth
\begin{equation*}
\widehat{h}_{\text{opt}}(h_1, q_1, q) = \widehat{c}_0(h_1, q_1, q)n^{1/(1+2q)}
\end{equation*}
in the definition of $\widehat{C}_{h,q}$ in~\eqref{eq:1}.
\end{enumerate}

\subsection{Dynamic functional principal component decomposition}

From the long-run covariance $C(w,v)$, we apply functional principal decomposition to extract the functional principal components and their associated scores. With Karhunen-Lo\`{e}ve expansion, a stochastic process $\X$ can be expressed as
\begin{equation*}
\X_t(\omega) = \mu(\omega) + \sum^{\infty}_{j=1}\beta_{t,j}\phi_j(\omega),
\end{equation*}
where $\X_t^c(\omega) = \X_t(\omega) - \mu(\omega)$ and $\beta_{t,j}$ is an uncorrelated random variable with zero mean and unit variance. The principal component score $\beta_{t,j}$ is given by the projection of $\X_t^c$ in the direction of the $j$\textsuperscript{th} eigenfunction $\phi_j$, i.e., $\beta_{t,j} = \langle \X_t^c(\omega), \phi_j(\omega)\rangle$. The scores $\bm{\beta}_j = (\beta_{1,j},\dots,\beta_{n,j})$ constitute an uncorrelated sequence of random variables with zero mean and variance $\lambda_j$ which is the $j$\textsuperscript{th} eigenvalue. They can be interpreted as the weights of the contribution of the functional principal components $\phi_j(\omega)$ to $\X_t^c(\omega)$.

Since the long-run covariance $C(w,v)$ is unknown, the population eigenvalues and eigenfunctions can only be approximated through realizations of $\X(\omega)$. The sample mean and sample covariance are given by
\begin{align*}
\overline{\X}(\omega) &= \frac{1}{n}\sum^n_{t=1}\X_t(\omega), \\
\widehat{C}(\omega, \nu) &= \sum^{\infty}_{j=1}\widehat{\lambda}_j\widehat{\phi}_j(\omega)\widehat{\phi}_j(\nu),
\end{align*}
where $\widehat{\lambda}_1>\widehat{\lambda}_2>\cdots \geq 0$ are the sample eigenvalues of $\widehat{C}(\omega, \nu)$, and $\big[\widehat{\phi}_1(\omega), \widehat{\phi}_2(\omega), \dots\big]$ are the corresponding orthogonal sample eigenfunctions. The realizations of the stochastic process $\X$ can be written as
\begin{equation*}
\X_t(\omega) = \widehat{\mu}(\omega)+\sum^{\infty}_{j=1}\widehat{\beta}_{t,j}\widehat{\phi}_j(\omega), \qquad t = 1,2,\dots,n, 
\end{equation*}
where $\widehat{\mu}(\omega) = \frac{1}{n}\sum^n_{t=1}\X_t(\omega)$, and $\widehat{\beta}_{t,j}$ is the $j$\textsuperscript{th} estimated principal component score for the $t$\textsuperscript{th} time period.

\cite{HKH15} showed that kernel sandwich estimator in~\eqref{eq:1} is a consistent estimator of the true and unknown long-run covariance, and estimated functional principal components and principal component scores extracted from the estimated long-run covariance are also consistent.  

\section{Hurst exponent estimators}\label{sec:3}

Let the first set of estimated dynamic principal component scores be $\bm{\widehat{\beta}}_1 = \{\widehat{\beta}_{1,1}, \widehat{\beta}_{2,1}, \dots,\widehat{\beta}_{n,1}\}$. Since we consider the first set of scores, we shall replace $\bm{\widehat{\beta}}_1$ by $\bm{\beta}$ hereafter. In \cite{LRS18}, they also consider a $L_2$ norm of multiple sets of scores and find the estimation results remain similar. Due to space constraints, we present our results based on the first set of principal component scores. With the univariate time series of scores $\bm{\beta}$, we evaluate and compare some Hurst exponent estimators from long-memory univariate time-series literature. 

The Hurst exponent can be estimated either via time- or frequency-domain based estimators. These estimators can be divided into parametric and semi-parametric ones. The theory of parametric estimators was developed by \cite{FT86} and \cite{Dahlhaus89}. Semiparametric estimators of the memory parameter have become popular since they do not require knowing the specific form of the short-memory structure. They are based on the periodograms of the series, and can be categorized into two types: the log-periodogram estimator first proposed by \cite{GPH83} and the local-Whittle estimator which is credited to \cite{Kuensch87} and further developed by \cite{Robinson95}. The log-periodogram estimator is akin to the ordinary least squares and the local-Whittle estimator to the maximum likelihood estimator in the frequency domain.

\subsection{Time-domain based estimators}

In Sections~\ref{sec:3.1} to~\ref{sec:3.5}, we present five methods based on a simple linear regression model. In Sections~\ref{sec:3.6} and~\ref{sec:3.7}, we present two methods based on the R/S estimator.

\subsubsection{Aggregated variance estimator}\label{sec:3.1}

The aggregated variance estimator is based on the property of self-similar processes that variances of the aggregated processes decrease at the rate $m^{2H-2}$ as the block size $m$ increases \citep[e.g.,][Section 4.4]{TTW95, TT97, Beran94}. Recall that for a long-range dependent linear process,
\begin{equation*}
\text{Var}(\overline{\beta})\sim C m^{2H-2},
\end{equation*}
where $C$ is a constant. Consequently,
\begin{equation*}
\log_{10}[\text{Var}(\overline{\beta})]\approx \log_{10}C + (2H-2)\log_{10}m
\end{equation*}
With the predictor variable of $\log_{10}m$ and the response variable of $\log_{10}[\text{Var}(\overline{\beta})]$, we apply a simple linear regression to obtain an estimate of the slope parameter. For instance, one may define the following procedure:
\begin{enumerate}
\item[1)] Divide the time series $\bm{\beta}$ into $K$ non-overlapping blocks with block size $m$ and then average within each block, that is considered the aggregated series
\begin{equation}
\overline{\beta}^{(m)}(k) = \frac{1}{m}\sum^{km}_{t=(k-1)m+1}\beta_{t}, \label{eq:aggvarfit_1}
\end{equation}
where $k=1,\dots,K$ denotes a block index and $K=n/m\geq 1$ denotes the number of blocks.
\item[2)] Compute the overall mean
\begin{equation*}
\overline{\beta}^{(m)} = \frac{1}{K}\sum^K_{k=1}\overline{\beta}^{(m)}(k)
\end{equation*}
\item[3)] For a given $m$, compute the sample variance of $\beta^{(m)}(k)$ as
\begin{equation}
\widehat{\text{Var}}(\overline{\beta}^{(m)}) = \frac{1}{K}\sum^K_{k=1}\left[\overline{\beta}^{(m)}(k)\right]^2 - \left(\overline{\beta}^{(m)}\right)^2. \label{eq:aggvarfit_2}
\end{equation}
\item[4)] Heuristically, when $m$ grows, $\text{Var}(\overline{\beta}^{(m)})\sim C m^{2H-2}$. Thus, $\text{Var}(\overline{\beta}^{(m)})$ grows approximately at the rate $m^{2H-2}$. For different values of $m=1,\dots,M$, compute~\eqref{eq:aggvarfit_1} and~\eqref{eq:aggvarfit_2} to obtain $\widehat{\text{Var}}(\overline{\bm{\beta}}) = \{\widehat{\text{Var}}\overline{\beta}^{(1)},\dots,\widehat{\text{Var}}\overline{\beta}^{(M)}\}$. It is recommended by \cite{TTW95} and \cite{TT97} to choose values of $m$ that are equispaced on a logarithmic scale. Then, regress $\log_{10}[\widehat{\text{Var}}(\bm{\beta})]$ against $\log_{10}(M)$ to obtain regression coefficient $\widehat{\theta}_{\text{aggvar}}$. The estimated value of $H$ is given by
\begin{equation*}
\widehat{H}_{\text{aggvar}} = \frac{\widehat{\theta}_{\text{aggvar}}+2}{2}.
\end{equation*}
\end{enumerate}

\subsubsection{Differencing variance estimator}\label{sec:3.2}

To distinguish non-stationarity from long-range dependence, we can difference the variance \citep[see, e.g.,][]{TT97}. For a given $m^*$, we compute the difference of the sample variance
\begin{equation}
\text{Var}(\overline{\beta}^{m^*}) = \widehat{\text{Var}}\overline{\beta}^{(m^*+1)} - \widehat{\text{Var}}\overline{\beta}^{(m^*)},\qquad m^*=1,\dots,M-1.\label{eq:difference_variance}
\end{equation}
For different values of $m^*$, compute~\eqref{eq:aggvarfit_1},~\eqref{eq:aggvarfit_2} and~\eqref{eq:difference_variance} to obtain $\widehat{\text{Var}}(\overline{\bm{\beta}}) = \{\widehat{\text{Var}}\overline{\beta}^{(1)},\dots, \widehat{\text{Var}}\overline{\beta}^{(M-1)}\}$. Then, regress $\log_{10}[\widehat{\text{Var}}(\overline{\bm{\beta}})]$ against $\log_{10}(M-1)$ to obtain regression coefficient $\widehat{\theta}_{\text{diffvar}}$. The estimated value of $H$ is given by
\begin{equation*}
\widehat{H}_{\text{diffvar}} = \frac{\widehat{\theta}_{\text{diffvar}}+2}{2}.
\end{equation*}

\subsubsection{Absolute values of the aggregated series}\label{sec:absval}

Similar to the aggregated variance, the data are split in the same fashion, and the aggregated mean is computed from~\eqref{eq:aggvarfit_1}. Instead of computing the sample variance, one finds the sum of the absolute values of the aggregated series, namely
\begin{equation}
\text{abs}(\overline{\beta}^{(m)}) = \frac{1}{K}\sum^K_{k=1}\left|\overline{\beta}^{(m)}(k)\right|. \label{eq:absval_1}
\end{equation}
For different values of $m=1,\dots,M$, compute~\eqref{eq:absval_1} to obtain $\text{abs}(\overline{\bm{\beta}}) = \{\text{abs}(\overline{\beta}^{(1)}), \dots, \text{abs}(\overline{\beta}^{(M)})\}$. Then, regress $\log_{10}[\text{abs}(\overline{\bm{\beta}})]$ against $\log_{10}(M)$ to obtain regression coefficient $\widehat{\theta}_{\text{absval}}$. The estimated value of $H$ is given by
\begin{equation*}
\widehat{H}_{\text{absval}} = \widehat{\theta}_{\text{absval}}+1.
\end{equation*}

\subsubsection{Higuchi's method}\label{sec:3.4}

Similar to the absolute values of the aggregated series, the method of \cite{Higuchi88} calculates the partial sums $Y(n)=\sum^n_{t=1}\beta_t$ of the time series $\bm{\beta}$, and then finding the normalized length of the curve, namely
\begin{equation}
L(m) = \frac{n-1}{m^3}\sum^m_{i=1}\frac{m}{n-i}\sum_{k=1}^{(n-i)/m}\left|Y(i+km)-Y[i+(k-1)m]\right|,\label{eq:higuchi}
\end{equation}
where $n$ is the sample size of the time series, $m$ is a block size and $[\cdot]$ denotes the greatest integer function. Since $\text{E}[L(m)]\sim Cm^{-D}$ where $D=2-H$. For different values of $m=1,\dots,M$, compute~\eqref{eq:higuchi} to obtain $\bm{L} = \{L(1),\dots, L(M)\}$. Then, regress $\log_{10}(\bm{L})$ against $\log_{10}(M)$ to obtain regression coefficient $\widehat{\theta}_{\text{Higuchi}}$. The estimated value of $H$ is given by
\begin{equation*}
\widehat{H}_{\text{Higuchi}} = \widehat{\theta}_{\text{Higuchi}}+2.
\end{equation*}

\subsubsection{Detrended fluctuation analysis (DFA)}\label{sec:3.5}

Also known as a variance of residuals or Peng's method, DFA was introduced by \cite{PBH+94} to provide evidence of long memory in deoxyribonucleic acid (DNA) sequences. It consists of the following steps:
\begin{enumerate}
\item[1)] The data series is divided into $K$ nonoverlapping blocks and each block with size $m$ such that $n=mK$. 
\item[2)] Within each of the $K$ blocks, we regress $T_l = \sum^l_{t=1}\beta_t$ against $l$ and estimate the variance of the residuals by
\begin{equation*}
S_m^2(k) = \frac{1}{m}\sum^{km}_{l=(k-1)m+1}(T_l - \widehat{\zeta}_{0,k}-\widehat{\zeta}_{1,k}l)^2,\qquad k=1,\dots,K,
\end{equation*}
where $\widehat{\zeta}_{0, k}$ and $\widehat{\zeta}_{1,k}$ are least squares regression estimates based on the $k$th block. 
\item[3)] Compute the average of the variance of the residuals
\begin{equation}
F^2(m) = \frac{1}{K}\sum^K_{k=1}S_m^2(k).\label{eq:Peng_method}
\end{equation}
\item[4)] Heuristically, $F^2(m)$ grows at the rate $m^{2H}$. For different values of $m=1,\dots,M$, compute~\eqref{eq:Peng_method} to obtain $\bm{F}^2 = \{F^2(1),\dots,F^2(M)\}$. Then, regress $\log_{10}\bm{F}^2$ against $\log_{10}M$ to obtain regression coefficient $\widehat{\theta}_{\text{Peng}}$. The estimated value of $H$ is given by
\begin{equation*}
\widehat{H}_{\text{Peng}} = \frac{\widehat{\theta}_{\text{Peng}}}{2}.
\end{equation*}
\end{enumerate}

The DFA bears a strong resemblance to the variance plot, but instead of assuming stationarity, a fitted linear trend is subtracted from each block \citep{BFG+13}. Therefore, the DFA is less sensitive to the trend exhibited in the data.

\subsubsection{Rescaled Range (R/S) estimator}\label{sec:3.6}

The R/S estimator was introduced by \cite{Hurst51} for estimating the minimum capacity of a dam. The R/S estimator is one of the first methods for estimating Hurst exponent. Although many Hurst exponent estimators have better statistical properties than the R/S estimator (which, for example, is inefficient in the case of Gaussian innovations), it is a simple method that computes fast \citep[see, e.g.,][]{LRS18}. Given a time series of scores $\bm{\beta}$, calculation of the R/S statistic has the following steps:
\begin{enumerate}
\item[1)] Calculate the range 
\begin{equation*}
R_n = \max_{1\leq T\leq n}\sum^{T}_{t=1}\left(\beta_{t} - \overline{\beta}\right) - \min_{1\leq T\leq n}\sum^{T}_{t=1}\left(\beta_{t} - \overline{\beta}\right), \qquad \overline{\beta} = \frac{1}{n}\sum^n_{t=1}\beta_{t}
\end{equation*}
\item[2)] Calculate the scale
\begin{equation*}
S_n = \sqrt{\frac{1}{n-1}\sum^n_{t=1}\left(\beta_{t} - \overline{\beta}\right)^2}.
\end{equation*}
If $\beta_t$ is second-order stationary, then $S_n^2$ converges in probability to $\text{Var}(\beta_t)$ \citep[][p.410]{BFG+13}.
\item[3)] The R/S estimator may be defined by
\begin{equation*}
\widehat{H}_{\text{RS}} = \frac{\log_{10}(R_n/S_n)}{\log_{10}n}.
\end{equation*}
The plot of $\log_{10}(R_n/S_n)$ against $\log_{10}n$ is also known as ``pox plots".
\end{enumerate}

\subsubsection{Rescaled adjusted range estimator}\label{sec:3.7}

While the R/S estimator is applied to the original time series, the rescaled adjusted range estimator is implemented to the partial sum of the original time series \citep[see, e.g.,][]{MW69, Mandelbrot75, MT79}. For a univariate time series of principal component scores $\bm{\beta}$ with the partial sum 
\begin{equation*}
Y(n) = \sum^n_{t=1}\beta_{t}
\end{equation*}
and sample variance 
\begin{equation*}
S^2(n):=\frac{1}{n}\sum^n_{t=1}\beta_{t}^2 - \left[\frac{1}{n}Y(n)\right]^2,
\end{equation*} 
the rescaled adjusted range estimator is given by
\begin{equation*}
\frac{R_n}{S_n} :=\frac{1}{S(n)}\left\{\max_{1\leq T\leq n}\left[Y(T)-\frac{T}{n}Y(n)\right]-\min_{1\leq T\leq n}\left[Y(T)-\frac{T}{n}Y(n)\right]\right\}.
\end{equation*}
Choosing logarithmically equidistant values of $n$, regress $\log_{10}(R_n/S_n)$ against $\log_{10}(n)$ to obtain regression coefficient $\widehat{\theta}_{\text{RAR}}$. The estimated value of $H$ is given by
\begin{equation*}
\widehat{H}_{\text{RAR}} = \widehat{\theta}_{\text{RAR}}.
\end{equation*}

\subsection{Frequency-domain based estimators}

\subsubsection{(Smoothed) periodogram estimator}

With a univariate time series of scores $\bm{\beta}$, the periodogram can be defined as 
\begin{equation}
w(\lambda_j) = (2\pi n)^{-1/2}\sum^n_{t=1}\beta_{t} \exp^{it\lambda_j}, \qquad I(\lambda_j) = |w(\lambda_j)|^2,\label{eq:period}
\end{equation}
where $\lambda_j = 2\pi j/n$ denotes the set of harmonic frequencies, $j=1,\dots,J$ where $J$ is a positive integer, and $i^2=-1$. Since the periodogram is a measure of autocovariance, it can also be expressed as
\begin{align*}
I(\lambda_j) &= \frac{1}{2\pi}\sum^{n-1}_{s=1-n}R(s)\cos(s\lambda_j) \\
&= \frac{1}{2\pi}\left\{R(0) +  2\sum^{n-1}_{s=1}R(s)\cos(s\lambda_j)\right\},\qquad \lambda_j\in [-\pi, \pi]
\end{align*}
where $R(s)$ denotes the sample autocovariance function, i.e.,
\begin{equation*}
R(s) = \frac{1}{n}\sum^{n-s}_{t=1}(\beta_t - \overline{\beta})(\beta_{t+s}-\overline{\beta}), \qquad s=0,\pm 1, \dots, \pm (n-1),
\end{equation*}
where $\overline{\beta}$ is the sample mean of the time series of scores.

Because $I(\lambda)$ is an estimator of the spectral density, a time series with long-range dependence should have a periodogram which is proportional to $|\lambda|^{1-2H} = |\lambda|^{-2d}$ close to the origin \citep{TTW95}. Thus, regress the logarithm of the periodogram for different values of $\lambda$ against $\log_{10}(\lambda)$ to obtain regression coefficient $\widehat{\theta}_{\text{per}}$. The estimated value of $H$ is given by
\begin{equation*}
\widehat{H}_{\text{per}} = \frac{1-\widehat{\theta}_{\text{per}}}{2}.
\end{equation*}
As advocated by \cite{TTW95}, we use only the lowest 10\% of the frequencies for the regression, since the proportionality above 10\% only holds for $\lambda$ close to the origin.

The frequency axis is divided into logarithmically equidistant boxes, and the periodogram values corresponding to the frequencies inside the box are averaged, to obtain smoothed periodogram. The periodogram values at very low frequencies are remained, while the rest are divided into 60 boxes \citep[see, e.g.,][]{TTW95}. By regressing the logarithm of the smoothed periodogram against frequencies, we obtain regression coefficient $\widehat{\theta}_{\text{boxper}}$. To achieve the robustness in the least square fitting, we use a robust linear model. The estimated value of $H$ is given by
\begin{equation*}
\widehat{H}_{\text{boxper}} = \frac{1-\widehat{\theta}_{\text{boxper}}}{2}.
\end{equation*}

\subsubsection{(Smoothed) Geweke-Porter-Hudak estimator}

In the univariate ARFIMA$(p, d, q)$ models, \cite{GPH83} proposed a semiparametric estimator of $d$ based on the first $J$ periodogram ordinates given in~\eqref{eq:period}. Let $\bm{\beta}$ be a stationary time series with spectral density
\begin{align}
f(\lambda) &= \big|1-\exp^{-i\lambda}\big|^{-2d}f_*(\lambda) \notag\\
&\sim |\lambda|^{-2d}f_*(\lambda) \notag\\
&\sim C |\lambda|^{-2d} \label{eq:GPH}
\end{align}
as $\lambda\rightarrow 0$, where $-\frac{1}{2}<d<\frac{1}{2}$. Recall that the empirical estimate to the spectral density is the periodogram given in~\eqref{eq:GPH},
\begin{equation*}
\log_{10}f(\lambda) \sim \log_{10}C + d b(\lambda),
\end{equation*}
where $b(\lambda) = -2\log_{10}(\lambda)$. In practice, we replace $f(\lambda)$ by its empirical analogy $I(\lambda)$, thus
\begin{equation*}
\log_{10}I(\lambda) \sim \log_{10}C + d b(\lambda). 
\end{equation*}

By a simple linear regression, \cite{GPH83} suggested the least-square estimator
\begin{align*}
\widehat{d}_{\text{GPH}} &= \frac{\sum^J_{j=1}(b_j-\overline{b})\log_{10} I(\lambda_j)}{\sum^J_{j=1}(b_j - \overline{b})^2}, \\
&=\frac{-\frac{1}{2}\sum^J_{j=1}[\log_{10}(\lambda_j) - \overline{\log_{10}(\lambda_j)}]\log_{10}I(\lambda_j)}{\sum^J_{j=1}[\log_{10}(\lambda_j) - \overline{\log_{10}(\lambda_j)}]^2}
\end{align*}
where $b_j=-2\log_{10}(\lambda_j)$, $\overline{b} = \frac{1}{J}\sum^J_{j=1}b_j$ and $\overline{\log_{10}(\lambda_j)}=\frac{1}{J}\sum^J_{j=1}\log_{10}(\lambda_j)$. Note that $\lambda_j = 2\pi j/n$ for $j=1,\dots, J$ are the $J$ smallest Fourier frequencies. The number $J$ acts as a bandwidth parameter. Following \cite{GPH83}, we choose $J=\sqrt{n}$.

Further, \cite{Robinson95b} showed that this estimator is consistent and has a central limit theorem of the form 
\begin{equation*}
\sqrt{J}(\widehat{d}_{\text{GPH}} - d) \xrightarrow[d]{} N\Big(0, \frac{\pi^2}{24}\Big).
\end{equation*}

\cite{Reisen94} considered a smoothed periodogram using the Parzen lag window, for estimating the parameter $d$. Let $I_s(\lambda)$ denote a smoothed periodogram of the form
\begin{equation*}
I_s(\lambda) = \frac{1}{2\pi}\sum^h_{s=-h}K\left(\frac{s}{h}\right)R(s)\cos(s\times \lambda),\qquad \lambda\in[-\pi, \pi],
\end{equation*}
where $K(u)$ is called the lag window generator, a fixed continuous even function in the range $-1<u<1$, with $K(0) = 1$ and $K(-u) = K(u)$. The bandwidth parameter $h$ is a function of $n$, and it is customarily chosen as $n^{\frac{9}{10}}$. The Parzen lag window generator has the following form:
\[ K(u) = \left\{ \begin{array}{ll}
         1-6u^2+6|u|^3 & \mbox{$|u| \leq \frac{1}{2}$};\\
        2(1-|u|)^3 & \mbox{$-\frac{1}{2} < u \leq 1$};\\
        0 & \mbox{$|u|>1$}.\end{array} \right. \] 

The smoothed periodogram estimator can be written as
\begin{align*}
\widehat{d}_{\text{SGPH}} &= \frac{\sum^J_{j=1}(b_j-\overline{b})\log_{10} I_s(\lambda_j)}{\sum^J_{j=1}(b_j - \overline{b})^2}, \\
&=\frac{-\frac{1}{2}\sum^J_{j=1}[\log_{10}(\lambda_j) - \overline{\log_{10}(\lambda_j)}]\log_{10}I_s(\lambda_j)}{\sum^J_{j=1}[\log_{10}(\lambda_j) - \overline{\log_{10}(\lambda_j)}]^2}.
\end{align*}

\subsubsection{Wavelet estimator}

This estimator computes the discrete wavelet transform, and obtains the wavelet coefficient $w_{j,k}$ associated with a mean zero $I(d)$ process with $-\frac{1}{2}<d<\frac{1}{2}$. The wavelet coefficient $w_{j,k}$ as $j\rightarrow 0$ are distributed $N(0, \sigma^22^{-2jd})$, where $\sigma^2$ is a finite constant \citep{Jensen99}. The variance $\sigma^2 2^{-2jd}$ depends on the scaling parameter $j$ but is independent of the translation parameter $k$. We define
\begin{equation*}
R(j) = \sigma^2 2^{-2jd}
\end{equation*}
be the wavelet coefficient's variance at scale $j$. Taking the logarithm transformation of $R(j)$, we obtain
\begin{equation*}
\log_{10}R(j) = \log_{10}\sigma^2 - d\log_{10}2^{2j},
\end{equation*}
where $d$ can be estimated via ordinary least squares. Since $R(j)$ is a population quantity, we estimate it by the sample variance of the wavelet coefficients as
\begin{equation*}
\widehat{R}(j) = \frac{1}{2^j}\sum^{2^j-1}_{k=0}w_{j,k}^2.
\end{equation*}

\subsubsection{Local Whittle estimator}

The local Whittle estimator is a Gaussian semiparametric estimation method to estimate the Hurst exponent based on the periodogram. It is first introduced by \cite{Kuensch87} and later developed by \cite{Robinson95}, \cite{Velasco99} and subsequent authors. The local Whittle method does not require the specification of a parametric model for the data. It only relies on the specification of the shape of the spectral density of the time series $\bm{\beta}$. 

Note that the spectral density $f(\lambda)$ of a stationary time series is usually assumed to satisfy that
\begin{equation*}
f(\lambda) \sim G \lambda^{1-2H}=G\lambda^{-2d}, \qquad \text{as}\quad \lambda\rightarrow 0+,
\end{equation*}
where $0<G<\infty$, $0<H<1$ and $-\frac{1}{2}<d<\frac{1}{2}$. 

Define $Q(G, d)$ as the objective function
\begin{equation}
Q(G, d) = \frac{1}{m_{\diamond}}\sum^{m_\diamond}_{j=1}\left\{\ln(G\lambda_j^{-2d})+\frac{I(\lambda_j)}{G\lambda_j^{-2d}}\right\},\label{eq:local_whittle}
\end{equation}
where $\lambda_j=(2\pi j)/n, j=1,\dots,m_{\diamond}$, and $m_\diamond$ is a positive integer satisfying $m_\diamond<n/2$, $m_{\diamond}\rightarrow \infty$ and $m_{\diamond}=o(n)$ \citep[see, e.g.,][]{Robinson95}. As in \cite{Robinson95}, we define the estimates
\begin{equation*}
\large(\widehat{G}, \widehat{d}\large) = \argmin_{0<G<\infty,\; d\in \Theta}Q(G, d),
\end{equation*}
where the closed interval of admissible estimates of true value of the self-similarity measure $d_0$, $\Theta = [\nabla_1, \nabla_2]$, $\nabla_1$ and $\nabla_2$ are numbers picked such that $-\frac{1}{2}<\nabla_1<\nabla_2<\frac{1}{2}$ as defined in \cite{Robinson95}. Alternatively, we may obtain
\begin{equation*}
\widehat{d} = \argmin_{d\in \Theta}R(d)
\end{equation*}
where
\begin{equation*}
R(d) = \ln \widehat{G}(d) - \frac{2d}{m_{\diamond}}\sum^{m_{\diamond}}_{j=1}\ln \lambda_j, \qquad \widehat{G}(d) = \frac{1}{m_\diamond}\sum^{m_\diamond}_{j=1}\lambda_j^{2d}I(\lambda_j).
\end{equation*}
Further, \cite{Robinson95} showed that $\widehat{d}$ is a consistent estimator of $d_0$, and $\sqrt{m_\diamond}(\widehat{d} - d_0)\rightarrow N(0, \frac{1}{4})$ as $n\rightarrow \infty$.

\subsubsection{Local Whittle estimator with tapering}

\cite{Velasco99} showed that it is possible to estimate consistently the Hurst exponent of non-stationary processes using the local Whittle estimator by tapering the observations. Let the tapered periodogram of $\bm{\beta}$ be $I_p(\lambda_j)$, and define $Q_p(G, d)$ as the objective function
\begin{equation}
Q_p(G, d) = \frac{p}{m_\diamond}\sum^{m_\diamond}_{j}\left\{\ln(G\lambda_j^{-2d})+\frac{I_p(\lambda_j)}{G\lambda_j^{-2d}}\right\},\label{eq:local_whittle_taper}
\end{equation}
where all the summations run for $j=p, 2p,\dots,m_\diamond$, assuming $m_\diamond/p$ is integer. Define the closed interval of admissible estimate of $d_0$, $\Theta = [\nabla_1, \nabla_2]$, $\nabla_1$ and $\nabla_2$ are numbers picked such that $0<\nabla_1<\nabla_2<d^*$ and $p\geq d^*+\frac{1}{2}$ where $d^*$ is the maximum value of $d$ we can estimate with tapers of order $p$, and $d^*$ may lie in a region where $\bm{\beta}$ is non-stationary. When $p=1$,~\eqref{eq:local_whittle_taper} reduces to~\eqref{eq:local_whittle}.

As in \cite{Velasco99}, we define the estimates
\begin{equation*}
(\widehat{G}_p, \widehat{d}_p) = \argmin_{0<G<\infty, \; d\in \Theta}Q_p(G, d).
\end{equation*}
Alternatively, we may obtain
\begin{equation*}
\widehat{d}_p = \argmin_{d\in\Theta}R_p(d).
\end{equation*}
where 
\begin{equation*}
R_p(d) = \ln \widehat{G}_p(d) - 2d\frac{p}{m_\diamond}\sum^{m_\diamond}_{j}\ln \lambda_j, \qquad \widehat{G}_p(d) = \frac{p}{m_\diamond}\sum^{m_\diamond}_{j}\lambda_j^{2d}I_p(\lambda_j).
\end{equation*}
The tapered periodogram includes only frequencies $\lambda_j$ for $j=p, 2p, \dots, m_\diamond$. The periodogram for non-stationary processes is equivalent to the periodogram for stationary processes evaluated at these frequencies \citep{Velasco99}.

\subsubsection{Modified local Whittle estimator}

\cite{HP14} proposed a modified local Whittle estimator that has good properties under local contamination. These contaminations include processes whose spectral density functions dominate at low frequencies, such as random level shifts, deterministic level shifts and deterministic trends \citep{HP14}. The data generating process is given
\begin{equation*}
z_t = c+\beta_t+u_t,
\end{equation*}
where $\beta_t$ is a process with memory parameter $d\in [0,\frac{1}{2}]$ and $c$ is a constant. When $d=0$, $\beta_t$ is a short-memory process. The process $u_t$ is the low frequency contamination. For a given sample size $n$, we define the periodogram of process $z_t$ to be $I_z(\lambda_j)$ and $f_z(\lambda_j)=\text{E}[I_z(\lambda_j)]$. Since the periodogram of $u_t$ is of order $O_p(\lambda_j^{-2}/n)$, we add a term $(G_u\lambda_j^{-2}/n)$ to the spectral density function of $\bm{\beta}$ to govern the low frequency contamination. The modified spectral density function is $f_j=G_0\lambda_j^{-2d}+G_u\lambda_j^{-2}/n$. Let $\theta=G_u/G_0$ be the noise-to-signal ratio, the modified spectral density function is
\begin{equation*}
    f_j = G_0(\lambda_j^{-2d}+\theta\lambda_j^{-2}/n) = G_0\cdot g_j,
\end{equation*}
where
\begin{equation*}
    g_j = \lambda_j^{-2d}+\theta\lambda_j^{-2}/n.
\end{equation*}
The modified local Whittle estimator is 
\begin{equation*}
    (\widehat{d}_{\diamond}, \widehat{\theta}_{\diamond}) = \argmin_{d,\; \theta}J_{m_\diamond}(d, \theta),
\end{equation*}
where
\begin{equation*}
    J_{m_\diamond}=\ln \left\{\frac{1}{m_\diamond}\sum^{m_\diamond}_{j=1}\left[\ln (g_j) + \frac{I_z(\lambda_j)}{g_j}\right]\right\}
\end{equation*}

\subsubsection{Exact local Whittle estimator}

The local Whittle estimator is based on an approximation of $I_{\beta}(\lambda_j)\sim \lambda_j^{-2d}I_u(\lambda_j)$, where $\beta$ denotes the original time series and $u$ denotes the noise process. \cite{SP05} proposed an exact local Whittle estimator that uses a corrected discrete Fourier transform of $\bm{\beta}$ to approximate periodogram $I_u(\lambda_j)$. They consider the fractional process $\bm{\beta}$ generated by the model
\begin{align*}
(1-L)^{d_0}\beta_t &= u_tI\{t\geq 1\}, \qquad t=0,\pm 1, \dots,\\
\beta_t &= (1-L)^{-d_0}u_tI\{t\geq 1\}.
\end{align*}
The discrete Fourier transform of a time series $\beta_t$ evaluated at frequency $\lambda$ as
\begin{align*}
w_{\beta}(\lambda) &= (2\pi n)^{-\frac{1}{2}}\sum^n_{t=1}\beta_t \exp^{it\lambda}\\
I_{\beta}(\lambda) &= |w_{\beta}(\lambda)|^2 \\
I_u(\lambda) &= I_{\Delta^{d_0}\beta}(\lambda)
\end{align*}
Define $Q(G, d)$ as the objective function
\begin{equation*}
Q_{m_\diamond}(G, d) = \frac{1}{m_{\diamond}}\sum^{m_\diamond}_{j=1}\left\{\ln(G\lambda_j^{-2d})+\frac{I_{\Delta^{d}\beta}(\lambda_j)}{G\lambda_j^{-2d}}\right\},
\end{equation*}
where $I_{\Delta^{d}\beta}(\lambda_j)$ is the periodogram of 
\begin{equation*}
\Delta^d\beta_t = (1-L)^d\beta_t = \sum^t_{k=0}\frac{(-d)_k}{k!}\beta_{t-k}.
\end{equation*}
As in \cite{Robinson95}, we define the estimates
\begin{equation*}
(\widehat{G}, \widehat{d}) = \argmin_{0<G<\infty, \; d\in \Theta}Q(G, d),
\end{equation*}
where $\Theta = [\Delta_1, \Delta_2]$, $\Delta_1$ and $\Delta_2$ are numbers picked such that $-\infty<\Delta_1<\Delta_2<\infty$. Alternatively, we may obtain
\begin{equation*}
\widehat{d} = \argmin_{d\in \Theta}R(d)
\end{equation*}
where
\begin{equation*}
R(d) = \ln \widehat{G}(d) - \frac{2d}{m_{\diamond}}\sum^{m_{\diamond}}_{j=1}\ln \lambda_j, \qquad \widehat{G}(H) = \frac{1}{m_\diamond}\sum^{m_\diamond}_{j=1}I_{\Delta^d \beta}(\lambda_j).
\end{equation*}

\subsubsection{Two-step local Whittle estimator}

The exact local Whittle estimator is consistent and has the same asymptotic distribution as the local Whittle estimator for all values of $d$ if the admissible range is less than 9/2 and the mean of the process is known. In practice, the mean of the process needs to be estimated, and \cite{Shimotsu10} studied the effect of an unknown mean on the exact local Whittle estimation. If an unknown mean is replaced by the simple average, then the exact local Whittle estimator is consistent for $d\in (-\frac{1}{2}, 1)$ and asymptotically normal for $d\in (-\frac{1}{2}, \frac{3}{4})$. 

\cite{Shimotsu10} considered the following data generating process
\begin{equation*}
\beta_t = \mu_0+\beta_t^0, \qquad \beta_t^0 = (1-L)^{-d_0}u_tI\{t\geq 1\},
\end{equation*}
where $\mu_0 = \text{E}(\beta_t)$ is a fixed unknown quantity. An estimator for $\mu_0$ is the sample average $\overline{\beta} = \frac{1}{n}\sum^n_{t=1}\beta_t$. The resulting memory parameter estimator is given as
\begin{equation*}
\widehat{d} = \argmin_{d\in \Theta}R^{\circ}(d),
\end{equation*}
where 
\begin{equation*}
R^{\circ}(d) = \ln \widehat{G}^{\circ}(d) - \frac{2d}{m_\diamond}\sum^{m_\diamond}_{j=1}\ln \lambda_j, \qquad \widehat{G}^{\circ}(d) = \frac{1}{m_\diamond}\sum^{m_\diamond}_{j=1}I_{\Delta^d(\bm{\beta}-\widehat{\mu})}(\lambda_j),
\end{equation*}
where $I_{\Delta^d (\bm{\beta}-\widehat{\mu})}(\lambda_j)$ is the periodogram of $\Delta^d(\bm{\beta} - \widehat{\mu})$.

\section{Numerical studies}\label{sec:4}

Numerical studies of finite-sample performance are provided via simulation and empirical applications.

\subsection{The functional ARFIMA model}

We study the functional ARFIMA$(p, d, q)$ process defined by
\begin{equation}
\nabla^d \X_t(u) = Y_t(u), \quad \nabla = 1-B, \quad -1/2<d<1/2, \label{eq:FARIMA_1}
\end{equation}
and
\begin{equation}
Y_t(u) - \sum^p_{i=1}\int_{\mathcal{I}}\phi_i(u, v)Y_{t-i}(v)dv = \eta_t(u) + \sum^q_{i=1}\int_{\mathcal{I}}\psi_i(u, v)\eta_{t-i}(v)dv,\label{eq:FARIMA_2}
\end{equation}
where $B$ denotes the backshift operator, $\{\eta_t\}$ denotes the noise operator, and $\phi_i(u, v)$ and $\psi_i(u, v)$ are the kernels with the associated integral operators defined by $\int_{\mathcal{I}}\phi_i(u, v)x(v)dv$ and $\int_{\mathcal{I}}\psi_i(u, v)x(v)dv$, respectively, $x\in \mathcal{H}$, and such that $Y_t(u)$ is stationary with respect to $t$. Note that
\begin{align}
\X_t &= \nabla^{-d}Y_t = (1-B)^{-d}Y_t \notag\\
&= \sum^{\infty}_{i=0}\beta_{i,-d}B^iY_t \notag\\
&= \sum^{\infty}_{i=0}\beta_{i,-d}Y_{t-i},\label{eq:MA_representation}
\end{align}
where, by Stirling's formula,
\begin{align*}
\beta_{i,-d} &= \beta_i^*+\beta_i^{\diamond} \\
\beta_i^* &= \frac{1}{\Gamma(d)}i^{-1+d} \\
\beta_i^{\diamond} &= O(i^{-2+d}),
\end{align*}
where $\Gamma(\cdot)$ is the gamma function. For a stationary series $Y_t$, it can be simulated from the following moving average (MA)($\infty$) representation of the functional autoregressive moving average (ARMA)$(p, q)$ process:
\begin{equation}
Y_t = \sum^{\infty}_{i=0}\pi[A_i(\overline{\eta}_{t-i})],\label{eq:ARMA}
\end{equation}
where $\pi(u_1,\dots, u_p)=u_1$, and $A_i$ denotes the integral operator in Hilbert space \citep[see, e.g.,][]{KKW17, LRS18}.
By combining~\eqref{eq:ARMA} with~\eqref{eq:MA_representation}, we obtain
\begin{equation*}
\X_t = \sum^{\infty}_{i=0}\beta_{i,-d}\sum^{\infty}_{j=0}\pi[A_j(\overline{\eta}_{t-i-j})].
\end{equation*} 
In our simulation studies, we implement a truncation, i.e., we use the first $n$ samples of $\overline{\eta}$ where $n$ denotes the sample size as the burn-in, and keep the remaining $n+100$ samples of $\overline{\eta}$ as our simulated realizations.

When $d=0$, model~\eqref{eq:FARIMA_1} and~\eqref{eq:FARIMA_2} becomes the functional ARMA$(p, q)$ of \cite{KKW17}, while when $q=0$, it further reduces to the functional autoregressive (AR)($p$) model of \cite{Bosq00} and \cite{LXC16}, and when $p=0$, it reduces to the functional MA$(q)$ model of \cite{CLT16} and \cite{AK17}. 

For notational simplicity, we let $\phi_i = \phi_i(\cdot, \cdot)$ and $\psi_i = \psi_i(\cdot, \cdot)$. As in \cite{Bosq00} and \cite{KKW17}, for the functional ARMA($p, q$) process $\{Y_t\}$, we can write
\begin{equation*}
\overline{\bm{Y}}_t(u) = \int_{\mathcal{I}}\overline{\bm{\phi}}(u, v)\overline{\bm{Y}}_{t-1}(v)dv+\sum^q_{i=0}\int_{\mathcal{I}}\overline{\bm{\psi}}_i(u, v)\overline{\bm{\eta}}_{t-i}(v)dv,
\end{equation*}
where
\begin{align*}
\overline{\bm{Y}}_t(u) &= [Y_t(u),\cdots, Y_{t-p+1}(u)]^{\top}, \qquad \overline{\bm{\eta}}_t(u) = [\eta_t(u), 0, \cdots, 0]^{\top}, \\
\overline{\bm{\phi}}(\cdot,\cdot) &= \begin{bmatrix}
    \phi_{1} & \phi_{2} & \cdots  & \phi_{p} \\
    I & O & \cdots & O \\
    \vdots & \vdots & \vdots  & \vdots \\
    O & \cdots & I & O
\end{bmatrix}, \qquad
\overline{\bm{\psi}}(\cdot,\cdot) = \begin{bmatrix}
    \psi_{i} & O & \cdots  & O \\
    O & O & \cdots & O \\
    \vdots & \vdots & \vdots  & \vdots \\
    O & O & \cdots & O
\end{bmatrix},
\end{align*}
$\psi_0=I$, $I$ and $O$ denote the identity and zero operators, respectively. 

\subsection{Simulation study}

We consider generating the curve time series $\X_t$ through a functional ARFIMA$(p, d, q)$ model defined in~\eqref{eq:FARIMA_1} and~\eqref{eq:FARIMA_2}, where $\mathcal{I} = [0,1]$, $\{\eta_t, t\in \mathbb{Z}\}$ is a sequence of independent and identically distributed standard Brownian motions over $[0, 1]$, and the following two cases are covered:
\begin{align*}
\text{Case 1:}& \quad p = 1, \ q = 0, \; \phi_1(u, v) = 0.34\times \exp\{-(u^2+v^2)/2\}, \\
\text{Case 2:}& \quad p = 1, \ q = 1, \; \phi_1(u, v) = 0.34\times \exp\{-(u^2+v^2)/2\}, \; \psi_1(u, v) = \frac{3}{2}\min(u, v),
\end{align*}
where $d=0.05, 0.10, \dots, 0.40$ in both cases. The choice of the constants in the definitions of $\phi_1$ and $\psi_1$ ensures that both $\|\phi_1\|$ and $\|\psi_1\|$ are smaller than one \citep[c.f.,][]{RS17}, so the simulated curve time series are stationary and invertible. The sample sizes employed are $n=250, 500, 1000$ with 1000 replications.

For a given estimation method, we obtain $B=1000$ estimated values of self-similar parameter $d$, namely $\widehat{d}_{b}$ for $b=1,\dots,1000$. We compute their bias, variance ($\sigma^2$) and mean squared error (MSE), given below
\begin{align*}
\text{Bias} &= \frac{1}{B}\sum_{b=1}^B (d - \widehat{d}_b), \\
\sigma^2 &= \frac{1}{B-1}\left[\sum^B_{b=1}\left(\widehat{d}_b - \frac{1}{B}\sum^B_{b=1}\widehat{d}_b\right)^2\right], \\
\text{MSE} &= \frac{1}{B}\sum^B_{b=1}(\widehat{d}_b - d)^2.
\end{align*}
The bias measures the tendency to over or under-estimate the long-memory parameter. The variance measures how far a set of estimated long-memory parameters are spread out from their mean. The MSE provides some information on the estimation accuracy of the long-memory parameter.

Under the functional ARFIMA($1, d, 0$) model, in Table~\ref{tab:1}, we evaluate and compare the finite-sample bias, variance, and MSE with the long-run covariance estimated from \cite{LRS18}. The R/S estimator produces the smallest variance. The local Whittle estimator with tapering produces the smallest bias, while the Peng's estimator produces the smallest MSE. 

\begin{center}
\tabcolsep 0.065in
\begin{longtable}{@{}lllrrrrrrrrr@{}}
\caption{With sample size $n=250, 500, 1000$, we evaluate and compare the finite-sample bias, variance, and MSE of the estimation error $\Delta d$ under the functional ARFIMA($1,d,0$) model, where the long-run covariance was estimated from \cite{LRS18}.}\label{tab:1}
\\
\toprule
 $n$ & Statistic & Estimator & \multicolumn{8}{c}{$d$} & Overall \\
& & & 0.05 & 0.10 & 0.15 & 0.20 & 0.25 & 0.30 & 0.35 & 0.40 & \\
\endfirsthead
\toprule
$n$ & Statistic & Estimator & \multicolumn{8}{c}{$d$} & Overall \\
& & & 0.05 & 0.10 & 0.15 & 0.20 & 0.25 & 0.30 & 0.35 & 0.40 & \\
\midrule
\endhead
\midrule
\multicolumn{12}{r}{Continued on next page} 
\endfoot
\endlastfoot
\midrule
250 & Bias & $\widehat{\epsilon}_{\text{aggvar}}$ & -0.087 & -0.098 & -0.109 & -0.121 & -0.134 & -0.148 & -0.164 & -0.182 & -0.130 \\ 
& & $\widehat{\epsilon}_{\text{diffvar}}$ & 0.228 & 0.224 & 0.223 & 0.227 & 0.228 & 0.232 & 0.239 & 0.226 & 0.228 \\ 
& & $\widehat{\epsilon}_{\text{absval}}$ & 0.028 & 0.017 & \textBF{0.006} & \textBF{-0.007} & -0.020 & -0.034 & -0.050 & -0.069 & -0.016 \\ 
& & $\widehat{\epsilon}_{\text{Higuchi}}$ & -0.044 & -0.053 & -0.061 & -0.071 & -0.082 & -0.095 & -0.111 & -0.132 & -0.081 \\ 
& & $\widehat{\epsilon}_{\text{Peng}}$ & 0.049 & 0.039 & 0.030 & 0.021 & 0.014 & 0.007 & \textBF{0.001} & \textBF{-0.004} & 0.020 \\ 
 & & $\widehat{\epsilon}_{\text{RS}}$ & 0.046 & \textBF{0.016} & -0.014 & -0.044 & -0.074 & -0.105 & -0.137 & -0.170 & -0.060 \\ 
& & $\widehat{\epsilon}_{\text{RAR}}$ & 0.120 & 0.094 & 0.066 & 0.039 & 0.011 & -0.018 & -0.047 & -0.077 & 0.023 \\ 
& & $\widehat{\epsilon}_{\text{per}}$ & 0.037 & 0.049 & 0.061 & 0.073 & 0.085 & 0.099 & 0.114 & 0.128 & 0.081 \\ 
& & $\widehat{\epsilon}_{\text{boxper}}$ & 0.132 & 0.132 & 0.132 & 0.132 & 0.132 & 0.131 & 0.132 & 0.134 & 0.132 \\ 
& & $\widehat{\epsilon}_{\text{GPH}}$   & 0.027 & 0.030 & 0.033 & 0.037 & 0.041 & 0.045 & 0.051 & 0.055 & 0.040 \\ 
& & $\widehat{\epsilon}_{\text{SGPH}}$ & -0.032 & -0.033 & -0.034 & -0.034 & -0.031 & -0.028 & -0.023 & -0.018 & -0.029 \\ 
& & $\widehat{\epsilon}_{\text{Wavelet}}$ & 0.052 & 0.048 & 0.044 & 0.040 & 0.036 & 0.032 & 0.028 & 0.025 & 0.038 \\ 
& & $\widehat{\epsilon}_{\text{Local\_W}}$ & -0.018 & -0.018 & -0.017 & -0.016 & -0.014 & -0.011 & -0.009 & -0.010 & -0.014 \\ 
& & $\widehat{\epsilon}_{\text{Local\_W\_T}}$ & 0.022 & 0.020 & 0.017 & 0.013 & \textBF{0.007} & \textBF{-0.002} & -0.013 & -0.029 & \textBF{0.004} \\ 
& & $\widehat{\epsilon}_{\text{Hou\_Perron}}$ &\textBF{-0.015} & -0.017 & -0.019 & -0.022 & -0.025 & -0.027 & -0.031 & -0.036 & -0.024 \\ 
& & $\widehat{\epsilon}_{\text{ELW}}$ & 0.060 & 0.061 & 0.062 & 0.063 & 0.064 & 0.066 & 0.068 & 0.070 & 0.064 \\ 
& & $\widehat{\epsilon}_{\text{ELW2S}}$ & 0.040 & 0.040 & 0.039 & 0.046 & 0.049 & 0.049 & 0.046 & 0.045 & 0.044 \\
\cmidrule{2-12}
& Variance & $\widehat{\epsilon}_{\text{aggvar}}$ & 0.030 & 0.029 & 0.029 & 0.029 & 0.028 & 0.027 & 0.026 & 0.025 & 0.028 \\ 
& &  $\widehat{\epsilon}_{\text{diffvar}}$ & 0.042 & 0.041 & 0.043 & 0.046 & 0.049 & 0.048 & 0.052 & 0.055 & 0.047 \\ 
& & $\widehat{\epsilon}_{\text{absval}}$ & 0.030 & 0.030 & 0.029 & 0.029 & 0.028 & 0.027 & 0.026 & 0.024 & 0.028 \\ 
& & $\widehat{\epsilon}_{\text{Higuchi}}$ & 0.015 & 0.016 & 0.017 & 0.018 & 0.019 & 0.020 & 0.020 & 0.019 & 0.018 \\ 
& & $\widehat{\epsilon}_{\text{Peng}}$ & 0.005 & 0.006 & 0.006 & 0.007 & 0.007 & 0.007 & 0.008 & 0.008 & 0.007 \\ 
 & & $\widehat{\epsilon}_{\text{RS}}$  & \textBF{0.002} & \textBF{0.002} & \textBF{0.002} & \textBF{0.002} & \textBF{0.002} & \textBF{0.002} & \textBF{0.002} & \textBF{0.002} & \textBF{0.002} \\ 
& & $\widehat{\epsilon}_{\text{RAR}}$ & 0.014 & 0.014 & 0.014 & 0.015 & 0.015 & 0.014 & 0.014 & 0.014 & 0.014 \\ 
& & $\widehat{\epsilon}_{\text{per}}$ & 0.036 & 0.035 & 0.034 & 0.033 & 0.033 & 0.032 & 0.032 & 0.032 & 0.033 \\ 
& & $\widehat{\epsilon}_{\text{boxper}}$ & 0.009 & 0.009 & 0.009 & 0.009 & 0.009 & 0.009 & 0.009 & 0.009 & 0.009 \\ 
& & $\widehat{\epsilon}_{\text{GPH}}$ & 0.045 & 0.044 & 0.043 & 0.042 & 0.042 & 0.041 & 0.041 & 0.041 & 0.042 \\ 
& & $\widehat{\epsilon}_{\text{SGPH}}$ & 0.024 & 0.025 & 0.025 & 0.025 & 0.025 & 0.025 & 0.026 & 0.026 & 0.025 \\ 
& & $\widehat{\epsilon}_{\text{Wavelet}}$ & 0.076 & 0.076 & 0.076 & 0.076 & 0.076 & 0.077 & 0.078 & 0.079 & 0.077 \\ 
& & $\widehat{\epsilon}_{\text{Local\_W}}$ & 0.009 & 0.009 & 0.009 & 0.009 & 0.009 & 0.009 & 0.008 & 0.007 & 0.009 \\ 
& & $\widehat{\epsilon}_{\text{Local\_W\_T}}$ & 0.045 & 0.044 & 0.042 & 0.039 & 0.036 & 0.031 & 0.027 & 0.022 & 0.036 \\ 
& & $\widehat{\epsilon}_{\text{Hou\_Perron}}$ & 0.015 & 0.016 & 0.016 & 0.018 & 0.019 & 0.021 & 0.023 & 0.027 & 0.019 \\
& & $\widehat{\epsilon}_{\text{ELW}}$ & 0.007 & 0.007 & 0.007 & 0.007 & 0.007 & 0.007 & 0.007 & 0.007 & 0.007 \\ 
& & $\widehat{\epsilon}_{\text{ELW2S}}$ & 0.020 & 0.020 & 0.019 & 0.020 & 0.019 & 0.018 & 0.018 & 0.020 & 0.019 \\
 \cmidrule{2-12}
& MSE & $\widehat{\epsilon}_{\text{aggvar}}$ & 0.037 & 0.039 & 0.041 & 0.043 & 0.046 & 0.049 & 0.053 & 0.058 & 0.046 \\ 
& &  $\widehat{\epsilon}_{\text{diffvar}}$ & 0.094 & 0.091 & 0.093 & 0.097 & 0.101 & 0.102 & 0.109 & 0.106 & 0.099 \\ 
& & $\widehat{\epsilon}_{\text{absval}}$ & 0.031 & 0.030 & 0.029 & 0.029 & 0.029 & 0.028 & 0.028 & 0.029 & 0.029 \\ 
& & $\widehat{\epsilon}_{\text{Higuchi}}$ & 0.017 & 0.019 & 0.021 & 0.023 & 0.026 & 0.029 & 0.032 & 0.037 & 0.025 \\ 
& & $\widehat{\epsilon}_{\text{Peng}}$ & 0.008 & 0.007 & 0.007 & 0.007 & \textBF{0.007} & \textBF{0.007} & \textBF{0.008} & 0.008 & \textBF{0.007} \\ 
 & & $\widehat{\epsilon}_{\text{RS}}$ & \textBF{0.004} & \textBF{0.002} & \textBF{0.002} & \textBF{0.004} & 0.007 & 0.013 & 0.021 & 0.031 & 0.011 \\ 
& & $\widehat{\epsilon}_{\text{RAR}}$ & 0.028 & 0.023 & 0.019 & 0.016 & 0.015 & 0.015 & 0.016 & 0.020 & 0.019 \\ 
& & $\widehat{\epsilon}_{\text{per}}$ & 0.037 & 0.037 & 0.037 & 0.038 & 0.040 & 0.042 & 0.045 & 0.049 & 0.041 \\ 
& & $\widehat{\epsilon}_{\text{boxper}}$ & 0.027 & 0.027 & 0.027 & 0.026 & 0.027 & 0.027 & 0.027 & 0.027 & 0.027 \\ 
& & $\widehat{\epsilon}_{\text{GPH}}$ & 0.046 & 0.045 & 0.044 & 0.043 & 0.044 & 0.043 & 0.044 & 0.044 & 0.044 \\ 
& & $\widehat{\epsilon}_{\text{SGPH}}$ & 0.025 & 0.026 & 0.026 & 0.026 & 0.026 & 0.026 & 0.026 & 0.027 & 0.026 \\ 
& & $\widehat{\epsilon}_{\text{Wavelet}}$ & 0.078 & 0.078 & 0.078 & 0.077 & 0.078 & 0.078 & 0.079 & 0.079 & 0.078 \\ 
& & $\widehat{\epsilon}_{\text{Local\_W}}$ & 0.009 & 0.009 & 0.009 & 0.009 & 0.009 & 0.009 & 0.009 & \textBF{0.007} & 0.009 \\ 
& & $\widehat{\epsilon}_{\text{Local\_W\_T}}$ & 0.049 & 0.048 & 0.045 & 0.042 & 0.037 & 0.032 & 0.027 & 0.022 & 0.038 \\ 
& & $\widehat{\epsilon}_{\text{Hou\_Perron}}$ & 0.015 & 0.016 & 0.017 & 0.018 & 0.019 & 0.021 & 0.024 & 0.028 & 0.020 \\
& & $\widehat{\epsilon}_{\text{ELW}}$ & 0.011 & 0.011 & 0.011 & 0.011 & 0.011 & 0.011 & 0.011 & 0.012 & 0.011 \\ 
& & $\widehat{\epsilon}_{\text{ELW2S}}$ & 0.021 & 0.021 & 0.021 & 0.022 & 0.021 & 0.021 &0.021 & 0.022 & 0.021 \\
\cmidrule{2-12}
500 & Bias & $\widehat{\epsilon}_{\text{aggvar}}$ & -0.068 & -0.078 & -0.090 & -0.102 & -0.115 & -0.130 & -0.146 & -0.162 & -0.111 \\ 
& &  $\widehat{\epsilon}_{\text{diffvar}}$ & 0.218 & 0.214 & 0.212 & 0.213 & 0.210 & 0.209 & 0.207 & 0.215 & 0.212 \\ 
& & $\widehat{\epsilon}_{\text{absval}}$ & 0.016 & 0.005 & -0.006 & -0.018 & -0.032 & -0.046 & -0.063 & -0.080 & -0.028 \\ 
& & $\widehat{\epsilon}_{\text{Higuchi}}$ & -0.036 & -0.043 & -0.051 & -0.059 & -0.069 & -0.082 & -0.097 & -0.115 & -0.069 \\ 
& & $\widehat{\epsilon}_{\text{Peng}}$ & 0.038 & 0.029 & 0.022 & 0.015 & 0.009 & \textBF{0.003} & \textBF{-0.001} & \textBF{-0.005} & 0.014 \\ 
 & & $\widehat{\epsilon}_{\text{RS}}$ & 0.042 & 0.014 & -0.013 & -0.041 & -0.070 & -0.099 & -0.129 & -0.160 & -0.057 \\ 
& & $\widehat{\epsilon}_{\text{RAR}}$ & 0.097 & 0.074 & 0.050 & 0.027 & 0.002 & -0.023 & -0.049 & -0.077 & 0.013 \\ 
& & $\widehat{\epsilon}_{\text{per}}$ & 0.022 & 0.028 & 0.036 & 0.044 & 0.051 & 0.059 & 0.068 & 0.078 & 0.048 \\ 
& & $\widehat{\epsilon}_{\text{boxper}}$ & 0.078 & 0.079 & 0.079 & 0.079 & 0.080 & 0.081 & 0.082 & 0.084 & 0.080 \\ 
& & $\widehat{\epsilon}_{\text{GPH}}$ & 0.012 & 0.013 & 0.016 & 0.019 & 0.021 & 0.025 & 0.030 & 0.035 & 0.021 \\ 
& & $\widehat{\epsilon}_{\text{SGPH}}$ & -0.027 & -0.029 & -0.029 & -0.029 & -0.027 & -0.025 & -0.022 & -0.017 & -0.026 \\ 
& & $\widehat{\epsilon}_{\text{Wavelet}}$ & 0.029 & 0.025 & 0.021 & 0.017 & 0.014 & 0.011 & 0.009 & 0.006 & 0.017 \\ 
& & $\widehat{\epsilon}_{\text{Local\_W}}$ & -0.023 & -0.023 & -0.022 & -0.022 & -0.020 & -0.019 & -0.017 & -0.016 & -0.020 \\ 
& & $\widehat{\epsilon}_{\text{Local\_W\_T}}$ & \textBF{0.003} & \textBF{0.002} & \textBF{0.002} & \textBF{0.001} & \textBF{-0.001} & -0.004 & -0.009 & -0.019 & \textBF{-0.003} \\ 
& & $\widehat{\epsilon}_{\text{Hou\_Perron}}$ & -0.020 & -0.021 & -0.022 & -0.024 & -0.027 & -0.029 & -0.032 & -0.034 & -0.026 \\  
& & $\widehat{\epsilon}_{\text{ELW}}$ & 0.032 & 0.032 & 0.033 & 0.033 & 0.034 & 0.035 & 0.036 & 0.037 & 0.034 \\ 
& & $\widehat{\epsilon}_{\text{ELW2S}}$ & 0.023 & 0.019 & 0.018 & 0.023 & 0.031 & 0.034 & 0.029 & 0.025 & 0.025   \\
\cmidrule{2-12}
& Variance & $\widehat{\epsilon}_{\text{aggvar}}$ & 0.016 & 0.017 & 0.017 & 0.018 & 0.018 & 0.017 & 0.017 & 0.016 & 0.017 \\ 
& &  $\widehat{\epsilon}_{\text{diffvar}}$ & 0.021 & 0.022 & 0.023 & 0.022 & 0.023 & 0.023 & 0.023 & 0.024 & 0.023 \\ 
& & $\widehat{\epsilon}_{\text{absval}}$ & 0.017 & 0.017 & 0.018 & 0.018 & 0.018 & 0.018 & 0.017 & 0.016 & 0.017 \\ 
& & $\widehat{\epsilon}_{\text{Higuchi}}$ & 0.010 & 0.011 & 0.012 & 0.013 & 0.014 & 0.014 & 0.014 & 0.013 & 0.013 \\ 
& & $\widehat{\epsilon}_{\text{Peng}}$ & 0.003 & 0.003 & 0.003 & 0.004 & 0.004 & 0.004 & 0.004 & 0.004 & 0.004 \\ 
 & & $\widehat{\epsilon}_{\text{RS}}$ & \textBF{0.001} & \textBF{0.001} & \textBF{0.002} & \textBF{0.002} & \textBF{0.002} & \textBF{0.002} & \textBF{0.002} & \textBF{0.002} & \textBF{0.002} \\ 
& & $\widehat{\epsilon}_{\text{RAR}}$ & 0.009 & 0.009 & 0.009 & 0.009 & 0.009 & 0.009 & 0.009 & 0.009 & 0.009 \\ 
& & $\widehat{\epsilon}_{\text{per}}$ & 0.014 & 0.014 & 0.014 & 0.014 & 0.014 & 0.014 & 0.013 & 0.013 & 0.014 \\ 
& & $\widehat{\epsilon}_{\text{boxper}}$ & 0.003 & 0.003 & 0.003 & 0.003 & 0.003 & 0.003 & 0.003 & 0.003 & 0.003 \\ 
& & $\widehat{\epsilon}_{\text{GPH}}$ & 0.028 & 0.028 & 0.027 & 0.027 & 0.027 & 0.027 & 0.027 & 0.027 & 0.027 \\ 
& & $\widehat{\epsilon}_{\text{SGPH}}$ & 0.016 & 0.016 & 0.016 & 0.016 & 0.017 & 0.017 & 0.017 & 0.017 & 0.017 \\ 
& & $\widehat{\epsilon}_{\text{Wavelet}}$ & 0.027 & 0.027 & 0.027 & 0.027 & 0.027 & 0.026 & 0.026 & 0.026 & 0.027 \\ 
& & $\widehat{\epsilon}_{\text{Local\_W}}$ & 0.006 & 0.006 & 0.006 & 0.006 & 0.006 & 0.006 & 0.005 & 0.005 & 0.005 \\ 
& & $\widehat{\epsilon}_{\text{Local\_W\_T}}$ & 0.023 & 0.023 & 0.023 & 0.022 & 0.021 & 0.020 & 0.018 & 0.015 & 0.020 \\ 
& & $\widehat{\epsilon}_{\text{Hou\_Perron}}$ & 0.008 & 0.008 & 0.009 & 0.009 & 0.011 & 0.011 & 0.013 & 0.014 & 0.010 \\  
& & $\widehat{\epsilon}_{\text{ELW}}$ & 0.004 & 0.004 & 0.004 & 0.004 & 0.004 & 0.004 & 0.004 & 0.004 & 0.004 \\ 
& & $\widehat{\epsilon}_{\text{ELW2S}}$ & 0.006 & 0.006 & 0.008 & 0.009 & 0.009 & 0.007 & 0.006 & 0.007 & 0.007 \\
\cmidrule{2-12}
& MSE & $\widehat{\epsilon}_{\text{aggvar}}$ & 0.021 & 0.023 & 0.025 & 0.028 & 0.031 & 0.034 & 0.038 & 0.042 & 0.030 \\ 
& &  $\widehat{\epsilon}_{\text{diffvar}}$ & 0.069 & 0.068 & 0.068 & 0.067 & 0.067 & 0.067 & 0.066 & 0.070 & 0.068 \\ 
& & $\widehat{\epsilon}_{\text{absval}}$ & 0.017 & 0.017 & 0.018 & 0.018 & 0.019 & 0.020 & 0.021 & 0.023 & 0.019 \\ 
& & $\widehat{\epsilon}_{\text{Higuchi}}$ & 0.011 & 0.013 & 0.015 & 0.017 & 0.019 & 0.021 & 0.023 & 0.027 & 0.018 \\ 
& & $\widehat{\epsilon}_{\text{Peng}}$ & 0.004 & 0.004 & 0.004 & 0.004 & \textBF{0.004} & \textBF{0.004} & \textBF{0.004} & \textBF{0.005} & \textBF{0.004} \\ 
 & & $\widehat{\epsilon}_{\text{RS}}$ & \textBF{0.003} & \textBF{0.002} & \textBF{0.002} & \textBF{0.003} & 0.007 & 0.012 & 0.018 & 0.027 & 0.009 \\ 
& & $\widehat{\epsilon}_{\text{RAR}}$ & 0.018 & 0.014 & 0.012 & 0.010 & 0.009 & 0.010 & 0.011 & 0.015 & 0.012 \\ 
& & $\widehat{\epsilon}_{\text{per}}$ & 0.015 & 0.015 & 0.015 & 0.016 & 0.016 & 0.017 & 0.018 & 0.020 & 0.016 \\ 
& & $\widehat{\epsilon}_{\text{boxper}}$ & 0.009 & 0.010 & 0.010 & 0.010 & 0.010 & 0.010 & 0.010 & 0.010 & 0.010 \\ 
& & $\widehat{\epsilon}_{\text{GPH}}$ & 0.028 & 0.028 & 0.027 & 0.027 & 0.028 & 0.028 & 0.028 & 0.028 & 0.028 \\ 
& & $\widehat{\epsilon}_{\text{SGPH}}$ & 0.017 & 0.017 & 0.017 & 0.017 & 0.017 & 0.017 & 0.017 & 0.017 & 0.017 \\ 
& & $\widehat{\epsilon}_{\text{Wavelet}}$ & 0.028 & 0.027 & 0.027 & 0.027 & 0.027 & 0.026 & 0.026 & 0.026 & 0.027 \\ 
& & $\widehat{\epsilon}_{\text{Local\_W}}$ & 0.006 & 0.006 & 0.006 & 0.006 & 0.006 & 0.006 & 0.006 & 0.005 & 0.006 \\ 
& & $\widehat{\epsilon}_{\text{Local\_W\_T}}$ & 0.023 & 0.023 & 0.023 & 0.022 & 0.021 & 0.020 & 0.018 & 0.015 & 0.020 \\ 
& & $\widehat{\epsilon}_{\text{Hou\_Perron}}$ & 0.009 & 0.009 & 0.009 & 0.010 & 0.011 & 0.012 & 0.014 & 0.015 & 0.011 \\
& & $\widehat{\epsilon}_{\text{ELW}}$ & 0.005 & 0.005 & 0.005 & 0.005 & 0.005 & 0.005 & 0.005 & 0.005 & 0.005 \\ 
& & $\widehat{\epsilon}_{\text{ELW2S}}$ & 0.007 & 0.007 & 0.008 & 0.010 & 0.010 & 0.008 & 0.007 & 0.007 & 0.008\\
\cmidrule{2-12}
1000 & Bias & $\widehat{\epsilon}_{\text{aggvar}}$ & -0.049 & -0.058 & -0.067 & -0.078 & -0.091 & -0.105 & -0.121 & -0.139 & -0.088 \\ 
& &  $\widehat{\epsilon}_{\text{diffvar}}$ & 0.191 & 0.188 & 0.189 & 0.189 & 0.188 & 0.187 & 0.188 & 0.188 & 0.188 \\ 
& & $\widehat{\epsilon}_{\text{absval}}$ & 0.014 & 0.005 & -0.004 & -0.015 & -0.028 & -0.042 & -0.058 & -0.076 & -0.026 \\ 
& & $\widehat{\epsilon}_{\text{Higuchi}}$ & -0.019 & -0.024 & -0.030 & -0.036 & -0.044 & -0.054 & -0.067 & -0.084 & -0.045 \\ 
& & $\widehat{\epsilon}_{\text{Peng}}$ & 0.030 & 0.023 & 0.017 & 0.011 & 0.007 & 0.002 & \textBF{-0.002} & -0.005 & 0.011 \\ 
 & & $\widehat{\epsilon}_{\text{RS}}$ & 0.039 & 0.014 & -0.011 & -0.037 & -0.063 & -0.090 & -0.119 & -0.148 & -0.052 \\ 
& & $\widehat{\epsilon}_{\text{RAR}}$ & 0.091 & 0.071 & 0.049 & 0.028 & 0.005 & -0.019 & -0.043 & -0.070 & 0.014 \\ 
& & $\widehat{\epsilon}_{\text{per}}$ & 0.017 & 0.021 & 0.026 & 0.031 & 0.036 & 0.041 & 0.047 & 0.053 & 0.034 \\ 
& & $\widehat{\epsilon}_{\text{boxper}}$ & 0.040 & 0.040 & 0.041 & 0.041 & 0.042 & 0.043 & 0.044 & 0.046 & 0.042 \\ 
& & $\widehat{\epsilon}_{\text{GPH}}$ & 0.007 & 0.008 & 0.010 & 0.012 & 0.015 & 0.017 & 0.021 & 0.026 & 0.014 \\ 
& & $\widehat{\epsilon}_{\text{SGPH}}$ & -0.021 & -0.022 & -0.022 & -0.022 & -0.020 & -0.018 & -0.015 & -0.011 & -0.019 \\ 
& & $\widehat{\epsilon}_{\text{Wavelet}}$ & 0.020 & 0.017 & 0.014 & 0.012 & 0.010 & 0.008 & 0.006 & \textBF{0.004} & 0.011 \\ 
& & $\widehat{\epsilon}_{\text{Local\_W}}$ & -0.013 & -0.013 & -0.012 & -0.012 & -0.010 & -0.009 & -0.007 & -0.005 & -0.010 \\ 
& & $\widehat{\epsilon}_{\text{Local\_W\_T}}$ & \textBF{-0.000} & \textBF{-0.000} & \textBF{-0.000} & \textBF{-0.000} & \textBF{-0.001} & \textBF{-0.002} & -0.005 & -0.011 & \textBF{-0.002} \\ 
& & $\widehat{\epsilon}_{\text{Hou\_Perron}}$ & -0.007 & -0.008 & -0.009 & -0.010 & -0.011 & -0.012 & -0.013 & -0.015 & -0.011 \\ 
& & $\widehat{\epsilon}_{\text{ELW}}$ & 0.022 & 0.022 & 0.022 & 0.023 & 0.023 & 0.024 & 0.025 & 0.026 & 0.023 \\ 
& & $\widehat{\epsilon}_{\text{ELW2S}}$ & 0.016 & 0.010 & 0.008 & 0.016 & 0.027 & 0.029 & 0.023 & 0.018 & 0.018 \\ 
\cmidrule{2-12}
& Variance & $\widehat{\epsilon}_{\text{aggvar}}$ & 0.011 & 0.012 & 0.012 & 0.012 & 0.012 & 0.012 & 0.011 & 0.011 & 0.012 \\ 
& &  $\widehat{\epsilon}_{\text{diffvar}}$ & 0.012 & 0.011 & 0.012 & 0.011 & 0.012 & 0.012 & 0.013 & 0.013 & 0.012 \\ 
& & $\widehat{\epsilon}_{\text{absval}}$ & 0.011 & 0.011 & 0.011 & 0.011 & 0.011 & 0.011 & 0.011 & 0.011 & 0.011 \\ 
& & $\widehat{\epsilon}_{\text{Higuchi}}$ & 0.005 & 0.005 & 0.006 & 0.006 & 0.007 & 0.007 & 0.007 & 0.007 & 0.006 \\ 
& & $\widehat{\epsilon}_{\text{Peng}}$ & 0.002 & 0.002 & 0.002 & 0.002 & 0.002 & 0.003 & 0.003 & 0.003 & 0.002 \\ 
 & & $\widehat{\epsilon}_{\text{RS}}$ & \textBF{0.001} & \textBF{0.001} & \textBF{0.001} & \textBF{0.001} & \textBF{0.002} & \textBF{0.002} & \textBF{0.002} & \textBF{0.002} & \textBF{0.001} \\ 
& & $\widehat{\epsilon}_{\text{RAR}}$ & 0.006 & 0.007 & 0.007 & 0.007 & 0.007 & 0.007 & 0.007 & 0.007 & 0.007 \\ 
& & $\widehat{\epsilon}_{\text{per}}$ & 0.006 & 0.006 & 0.006 & 0.006 & 0.006 & 0.006 & 0.006 & 0.006 & 0.006 \\ 
& & $\widehat{\epsilon}_{\text{boxper}}$ & 0.002 & 0.002 & 0.002 & 0.002 & 0.002 & 0.002 & 0.002 & 0.002 & 0.002 \\ 
& & $\widehat{\epsilon}_{\text{GPH}}$ & 0.019 & 0.019 & 0.019 & 0.019 & 0.019 & 0.019 & 0.019 & 0.020 & 0.019 \\ 
& & $\widehat{\epsilon}_{\text{SGPH}}$ & 0.011 & 0.011 & 0.011 & 0.011 & 0.012 & 0.012 & 0.012 & 0.012 & 0.012 \\ 
& & $\widehat{\epsilon}_{\text{Wavelet}}$ & 0.012 & 0.012 & 0.012 & 0.012 & 0.012 & 0.012 & 0.012 & 0.012 & 0.012 \\ 
& & $\widehat{\epsilon}_{\text{Local\_W}}$ & 0.003 & 0.003 & 0.003 & 0.003 & 0.003 & 0.003 & 0.003 & 0.003 & 0.003 \\ 
& & $\widehat{\epsilon}_{\text{Local\_W\_T}}$ & 0.014 & 0.014 & 0.014 & 0.014 & 0.013 & 0.013 & 0.012 & 0.010 & 0.013 \\ 
& & $\widehat{\epsilon}_{\text{Hou\_Perron}}$ & 0.004 & 0.004 & 0.005 & 0.005 & 0.005 & 0.005 & 0.006 & 0.006 & 0.005 \\ 
& & $\widehat{\epsilon}_{\text{ELW}}$ & 0.002 & 0.002 & 0.002 & 0.002 & 0.002 & 0.002 & 0.002 & 0.002 & 0.002 \\ 
& & $\widehat{\epsilon}_{\text{ELW2S}}$ & 0.003 & 0.003 & 0.004 & 0.006 & 0.006 & 0.005 & 0.004 & 0.004 & 0.004 \\ 
\cmidrule{2-12}
& MSE & $\widehat{\epsilon}_{\text{aggvar}}$ & 0.014 & 0.015 & 0.016 & 0.018 & 0.020 & 0.022 & 0.026 & 0.030 & 0.020 \\ 
& &  $\widehat{\epsilon}_{\text{diffvar}}$ & 0.048 & 0.046 & 0.047 & 0.047 & 0.047 & 0.047 & 0.048 & 0.049 & 0.047 \\ 
& & $\widehat{\epsilon}_{\text{absval}}$ & 0.012 & 0.012 & 0.012 & 0.012 & 0.013 & 0.014 & 0.015 & 0.017 & 0.013 \\ 
& & $\widehat{\epsilon}_{\text{Higuchi}}$ & 0.005 & 0.006 & 0.007 & 0.008 & 0.009 & 0.010 & 0.011 & 0.014 & 0.009 \\ 
& & $\widehat{\epsilon}_{\text{Peng}}$ & \textBF{0.003} & 0.002 & 0.002 & \textBF{0.002} & \textBF{0.002} & \textBF{0.003} & \textBF{0.003} & \textBF{0.003} & \textBF{0.003} \\ 
 & & $\widehat{\epsilon}_{\text{RS}}$ & 0.003 & \textBF{0.001} & \textBF{0.001} & 0.003 & 0.006 & 0.010 & 0.016 & 0.024 & 0.008 \\ 
& & $\widehat{\epsilon}_{\text{RAR}}$ & 0.015 & 0.012 & 0.009 & 0.008 & 0.007 & 0.008 & 0.009 & 0.012 & 0.010 \\ 
& & $\widehat{\epsilon}_{\text{per}}$ & 0.006 & 0.006 & 0.007 & 0.007 & 0.007 & 0.008 & 0.008 & 0.009 & 0.007 \\ 
& & $\widehat{\epsilon}_{\text{boxper}}$ & 0.003 & 0.003 & 0.003 & 0.003 & 0.003 & 0.003 & 0.004 & 0.004 & 0.003 \\ 
& & $\widehat{\epsilon}_{\text{GPH}}$ & 0.019 & 0.019 & 0.019 & 0.019 & 0.019 & 0.019 & 0.020 & 0.020 & 0.019 \\ 
& & $\widehat{\epsilon}_{\text{SGPH}}$ & 0.012 & 0.012 & 0.012 & 0.012 & 0.012 & 0.012 & 0.012 & 0.012 & 0.012 \\ 
& & $\widehat{\epsilon}_{\text{Wavelet}}$ & 0.012 & 0.012 & 0.012 & 0.012 & 0.012 & 0.012 & 0.012 & 0.012 & 0.012 \\ 
& & $\widehat{\epsilon}_{\text{Local\_W}}$ & 0.004 & 0.003 & 0.003 & 0.003 & 0.003 & 0.003 & 0.003 & 0.003 & 0.003 \\ 
& & $\widehat{\epsilon}_{\text{Local\_W\_T}}$ & 0.014 & 0.014 & 0.014 & 0.014 & 0.013 & 0.013 & 0.012 & 0.010 & 0.013 \\ 
& & $\widehat{\epsilon}_{\text{Hou\_Perron}}$ & 0.004 & 0.005 & 0.005 & 0.005 & 0.005 & 0.005 & 0.006 & 0.007 & 0.005 \\ 
& & $\widehat{\epsilon}_{\text{ELW}}$ & 0.003 & 0.003 & 0.003 & 0.003 & 0.003 & 0.003 & 0.003 & 0.003 & 0.003 \\ 
& & $\widehat{\epsilon}_{\text{ELW2S}}$ & 0.003 & 0.003 & 0.004 & 0.007 & 0.007 & 0.006 & 0.004 & 0.004 & 0.005 \\  
\bottomrule
\end{longtable}
\end{center}

Under the functional ARFIMA($1, d, 1$) model, in Table~\ref{tab:11}, we evaluate and compare the finite-sample bias, variance, and MSE with the long-run covariance estimated from \cite{LRS18}. The R/S estimator produces the smallest variance. The local Whittle and Hou-Perron estimators produce the smallest bias for various sample sizes, while the local Whittle estimator produces the smallest MSE.

\begin{center}
\tabcolsep 0.06in
\begin{longtable}{@{}lllrrrrrrrrr@{}}
\caption{With sample size $n=250, 500, 1000$, we evaluate and compare the finite-sample bias, variance, and MSE of the estimation error $\Delta d$ under the functional ARFIMA($1,d,1$) model, where the long-run covariance was estimated from \cite{LRS18}.}\label{tab:11}
\\
\toprule
 $n$ & Statistic & Estimator & \multicolumn{8}{c}{$d$} & Overall \\
& & & 0.05 & 0.10 & 0.15 & 0.20 & 0.25 & 0.30 & 0.35 & 0.40 & \\
\endfirsthead
\toprule
$n$ & Statistic & Estimator & \multicolumn{8}{c}{$d$} & Overall \\
& & & 0.05 & 0.10 & 0.15 & 0.20 & 0.25 & 0.30 & 0.35 & 0.40 & \\
\midrule
\endhead
\midrule
\multicolumn{12}{r}{Continued on next page} 
\endfoot
\endlastfoot
\midrule
250 & Bias & $\widehat{\epsilon}_{\text{aggvar}}$ & -0.091 & -0.104 & -0.117 & -0.131 & -0.146 & -0.162 & -0.179 & -0.198 & -0.141 \\ 
& &  $\widehat{\epsilon}_{\text{diffvar}}$ & 0.240 & 0.237 & 0.234 & 0.234 & 0.226 & 0.230 & 0.226 & 0.221 & 0.231 \\ 
& & $\widehat{\epsilon}_{\text{absval}}$ & 0.024 & 0.011 & \textBF{-0.002} & -0.017 & -0.032 & -0.049 & -0.066 & -0.086 & -0.027 \\ 
& & $\widehat{\epsilon}_{\text{Higuchi}}$ & -0.035 & -0.045 & -0.055 & -0.066 & -0.079 & -0.094 & -0.111 & -0.132 & -0.077 \\ 
& & $\widehat{\epsilon}_{\text{Peng}}$ & 0.106 & 0.094 & 0.083 & 0.072 & 0.062 & 0.053 & 0.045 & 0.038 & 0.069 \\ 
 & & $\widehat{\epsilon}_{\text{RS}}$ & 0.069 & 0.037 & 0.004 & -0.028 & -0.061 & -0.095 & -0.129 & -0.165 & -0.046 \\ 
& & $\widehat{\epsilon}_{\text{RAR}}$ & 0.146 & 0.117 & 0.087 & 0.056 & 0.025 & \textBF{-0.006} & -0.038 & -0.071 & 0.040 \\ 
& & $\widehat{\epsilon}_{\text{per}}$ & 0.046 & 0.057 & 0.068 & 0.080 & 0.091 & 0.104 & 0.118 & 0.132 & 0.087 \\ 
& & $\widehat{\epsilon}_{\text{boxper}}$ & 0.404 & 0.403 & 0.402 & 0.401 & 0.399 & 0.398 & 0.395 & 0.393 & 0.399 \\ 
& & $\widehat{\epsilon}_{\text{GPH}}$ & 0.022 & 0.025 & 0.027 & 0.029 & 0.031 & 0.034 & 0.039 & 0.043 & 0.031 \\ 
& & $\widehat{\epsilon}_{\text{SGPH}}$ & -0.033 & -0.036 & -0.037 & -0.038 & -0.037 & -0.035 & -0.031 & -0.027 & -0.034 \\ 
& & $\widehat{\epsilon}_{\text{Wavelet}}$ & 0.167 & 0.161 & 0.156 & 0.150 & 0.144 & 0.138 & 0.133 & 0.127 & 0.147 \\ 
& & $\widehat{\epsilon}_{\text{Local\_W}}$ & \textBF{0.007} & \textBF{0.007} & 0.008 & \textBF{0.009} & \textBF{0.010} & 0.011 & 0.012 & 0.008 & \textBF{0.009} \\ 
& & $\widehat{\epsilon}_{\text{Local\_W\_T}}$ & 0.060 & 0.057 & 0.053 & 0.047 & 0.039 & 0.028 & 0.013 & \textBF{-0.007} & 0.036 \\ 
& & $\widehat{\epsilon}_{\text{Hou\_Perron}}$ & 0.021 & 0.019 & 0.018 & 0.016 & 0.014 & 0.013 & \textBF{0.010} & 0.007 & 0.015 \\ 
& & $\widehat{\epsilon}_{\text{ELW}}$ & 0.110 & 0.110 & 0.111 & 0.112 & 0.113 & 0.114 & 0.115 & 0.117 & 0.113 \\ 
& & $\widehat{\epsilon}_{\text{ELW2S}}$ & 0.093 & 0.092 & 0.099 & 0.103 & 0.102 & 0.099 & 0.099 & 0.103 & 0.099 \\ 
\cmidrule{2-12}
& Variance & $\widehat{\epsilon}_{\text{aggvar}}$ & 0.031 & 0.031 & 0.031 & 0.031 & 0.031 & 0.031 & 0.030 & 0.028 & 0.031 \\ 
& &  $\widehat{\epsilon}_{\text{diffvar}}$ & 0.043 & 0.044 & 0.046 & 0.050 & 0.052 & 0.050 & 0.050 & 0.052 & 0.048 \\ 
& & $\widehat{\epsilon}_{\text{absval}}$ & 0.031 & 0.032 & 0.032 & 0.032 & 0.031 & 0.031 & 0.030 & 0.028 & 0.031 \\ 
& & $\widehat{\epsilon}_{\text{Higuchi}}$ & 0.016 & 0.017 & 0.018 & 0.019 & 0.020 & 0.020 & 0.020 & 0.020 & 0.019 \\ 
& & $\widehat{\epsilon}_{\text{Peng}}$ & 0.006 & 0.006 & 0.006 & 0.007 & 0.007 & 0.007 & 0.008 & 0.008 & 0.007 \\ 
 & & $\widehat{\epsilon}_{\text{RS}}$ & \textBF{0.002} & \textBF{0.002} & \textBF{0.002} & \textBF{0.002} & \textBF{0.002} & \textBF{0.002} & \textBF{0.002} & \textBF{0.002} & \textBF{0.002} \\ 
& & $\widehat{\epsilon}_{\text{RAR}}$ & 0.012 & 0.012 & 0.012 & 0.012 & 0.012 & 0.012 & 0.012 & 0.011 & 0.012 \\ 
& & $\widehat{\epsilon}_{\text{per}}$ & 0.037 & 0.036 & 0.036 & 0.036 & 0.035 & 0.035 & 0.035 & 0.035 & 0.036 \\ 
& & $\widehat{\epsilon}_{\text{boxper}}$ & 0.009 & 0.009 & 0.009 & 0.009 & 0.009 & 0.009 & 0.010 & 0.009 & 0.009 \\ 
& & $\widehat{\epsilon}_{\text{GPH}}$ & 0.047 & 0.047 & 0.046 & 0.046 & 0.046 & 0.045 & 0.045 & 0.046 & 0.046 \\ 
& & $\widehat{\epsilon}_{\text{SGPH}}$ & 0.025 & 0.026 & 0.026 & 0.027 & 0.027 & 0.028 & 0.028 & 0.029 & 0.027 \\ 
& & $\widehat{\epsilon}_{\text{Wavelet}}$ & 0.077 & 0.076 & 0.075 & 0.075 & 0.075 & 0.076 & 0.076 & 0.076 & 0.076 \\ 
& & $\widehat{\epsilon}_{\text{Local\_W}}$ & 0.009 & 0.009 & 0.009 & 0.009 & 0.009 & 0.009 & 0.008 & 0.007 & 0.009 \\ 
& & $\widehat{\epsilon}_{\text{Local\_W\_T}}$ & 0.045 & 0.044 & 0.042 & 0.039 & 0.036 & 0.031 & 0.027 & 0.022 & 0.036 \\ 
& & $\widehat{\epsilon}_{\text{Hou\_Perron}}$ & 0.014 & 0.014 & 0.015 & 0.016 & 0.017 & 0.018 & 0.020 & 0.022 & 0.017 \\ 
& & $\widehat{\epsilon}_{\text{ELW}}$ & 0.008 & 0.008 & 0.008 & 0.008 & 0.008 & 0.008 & 0.007 & 0.007 & 0.008 \\ 
& & $\widehat{\epsilon}_{\text{ELW2S}}$ & 0.019 & 0.019 & 0.020 & 0.019 & 0.018 & 0.018 & 0.019 & 0.020 & 0.019 \\
\cmidrule{2-12}
& MSE & $\widehat{\epsilon}_{\text{aggvar}}$ & 0.040 & 0.042 & 0.045 & 0.049 & 0.053 & 0.057 & 0.062 & 0.067 & 0.052 \\ 
& &  $\widehat{\epsilon}_{\text{diffvar}}$ & 0.101 & 0.100 & 0.101 & 0.104 & 0.102 & 0.103 & 0.101 & 0.101 & 0.102 \\ 
& & $\widehat{\epsilon}_{\text{absval}}$ & 0.032 & 0.032 & 0.032 & 0.032 & 0.032 & 0.033 & 0.034 & 0.035 & 0.033 \\ 
& & $\widehat{\epsilon}_{\text{Higuchi}}$ & 0.017 & 0.019 & 0.021 & 0.023 & 0.026 & 0.029 & 0.033 & 0.038 & 0.026 \\ 
& & $\widehat{\epsilon}_{\text{Peng}}$ & 0.017 & 0.015 & 0.013 & 0.012 & 0.011 & 0.010 & 0.010 & 0.010 & 0.012 \\ 
 & & $\widehat{\epsilon}_{\text{RS}}$ & \textBF{0.007} & \textBF{0.003} & \textBF{0.002} & \textBF{0.003} & \textBF{0.006} & 0.011 & 0.019 & 0.029 & 0.010 \\ 
& & $\widehat{\epsilon}_{\text{RAR}}$ & 0.034 & 0.026 & 0.020 & 0.016 & 0.013 & 0.012 & 0.013 & 0.016 & 0.019 \\ 
& & $\widehat{\epsilon}_{\text{per}}$ & 0.039 & 0.040 & 0.041 & 0.042 & 0.044 & 0.046 & 0.049 & 0.053 & 0.044 \\ 
& & $\widehat{\epsilon}_{\text{boxper}}$ & 0.173 & 0.172 & 0.171 & 0.170 & 0.169 & 0.168 & 0.166 & 0.164 & 0.169 \\ 
& & $\widehat{\epsilon}_{\text{GPH}}$ & 0.048 & 0.047 & 0.047 & 0.046 & 0.046 & 0.046 & 0.047 & 0.048 & 0.047 \\ 
& & $\widehat{\epsilon}_{\text{SGPH}}$ & 0.026 & 0.027 & 0.028 & 0.028 & 0.029 & 0.029 & 0.029 & 0.030 & 0.028 \\ 
& & $\widehat{\epsilon}_{\text{Wavelet}}$ & 0.104 & 0.102 & 0.099 & 0.097 & 0.096 & 0.095 & 0.093 & 0.092 & 0.097 \\ 
& & $\widehat{\epsilon}_{\text{Local\_W}}$ & 0.009 & 0.009 & 0.009 & 0.009 & 0.009 & \textBF{0.009} & \textBF{0.009} & \textBF{0.007} & \textBF{0.009} \\ 
& & $\widehat{\epsilon}_{\text{Local\_W\_T}}$ & 0.049 & 0.048 & 0.045 & 0.042 & 0.037 & 0.032 & 0.027 & 0.022 & 0.038 \\ 
& & $\widehat{\epsilon}_{\text{Hou\_Perron}}$ & 0.014 & 0.015 & 0.015 & 0.016 & 0.017 & 0.018 & 0.020 & 0.022 & 0.017 \\ 
& & $\widehat{\epsilon}_{\text{ELW}}$ & 0.020 & 0.020 & 0.020 & 0.020 & 0.020 & 0.020 & 0.021 & 0.021 & 0.020 \\ 
& & $\widehat{\epsilon}_{\text{ELW2S}}$ & 0.028 & 0.028 & 0.029 & 0.029 & 0.028 & 0.028 & 0.029 & 0.031 & 0.029 \\ 
\cmidrule{2-12}
500 & Bias & $\widehat{\epsilon}_{\text{aggvar}}$ & -0.068 & -0.080 & -0.092 & -0.106 & -0.120 & -0.136 & -0.153 & -0.170 & -0.115 \\ 
& &  $\widehat{\epsilon}_{\text{diffvar}}$ & 0.228 & 0.222 & 0.222 & 0.223 & 0.214 & 0.214 & 0.213 & 0.215 & 0.219 \\ 
& & $\widehat{\epsilon}_{\text{absval}}$ & 0.016 & 0.004 & -0.009 & -0.022 & -0.037 & -0.052 & -0.070 & -0.087 & -0.032 \\ 
& & $\widehat{\epsilon}_{\text{Higuchi}}$ & -0.028 & -0.037 & -0.046 & -0.055 & -0.066 & -0.079 & -0.095 & -0.114 & -0.065 \\ 
& & $\widehat{\epsilon}_{\text{Peng}}$ & 0.083 & 0.073 & 0.063 & 0.055 & 0.047 & 0.039 & 0.033 & 0.027 & 0.052 \\ 
 & & $\widehat{\epsilon}_{\text{RS}}$ & 0.064 & 0.034 & 0.004 & -0.026 & -0.057 & -0.089 & -0.121 & -0.154 & -0.043 \\ 
& & $\widehat{\epsilon}_{\text{RAR}}$ & 0.122 & 0.096 & 0.069 & 0.042 & 0.014 & -0.014 & -0.044 & -0.074 & 0.026 \\ 
& & $\widehat{\epsilon}_{\text{per}}$ & 0.033 & 0.039 & 0.046 & 0.054 & 0.061 & 0.069 & 0.078 & 0.089 & 0.059 \\ 
& & $\widehat{\epsilon}_{\text{boxper}}$ & 0.286 & 0.285 & 0.285 & 0.285 & 0.285 & 0.285 & 0.285 & 0.285 & 0.285 \\ 
& & $\widehat{\epsilon}_{\text{GPH}}$ & 0.010 & 0.011 & 0.013 & 0.014 & 0.017 & 0.020 & 0.024 & 0.029 & 0.017 \\ 
& & $\widehat{\epsilon}_{\text{SGPH}}$ & -0.029 & -0.031 & -0.032 & -0.032 & -0.031 & -0.029 & -0.026 & -0.021 & -0.029 \\ 
& & $\widehat{\epsilon}_{\text{Wavelet}}$ & 0.108 & 0.102 & 0.097 & 0.092 & 0.088 & 0.084 & 0.079 & 0.075 & 0.091 \\ 
& & $\widehat{\epsilon}_{\text{Local\_W}}$ & -0.009 & -0.009 & -0.009 & -0.008 & -0.007 & \textBF{-0.006} & \textBF{-0.005} & -0.005 & -0.007 \\ 
& & $\widehat{\epsilon}_{\text{Local\_W\_T}}$ & 0.022 & 0.022 & 0.021 & 0.020 & 0.018 & 0.014 & 0.007 & \textBF{-0.005} & 0.015 \\ 
& & $\widehat{\epsilon}_{\text{Hou\_Perron}}$ & \textBF{-0.000} & \textBF{-0.002} & \textBF{-0.002} & \textBF{-0.003} & \textBF{-0.005} & -0.007 & -0.010 & -0.011 & \textBF{-0.005} \\ 
& & $\widehat{\epsilon}_{\text{ELW}}$ & 0.062 & 0.062 & 0.062 & 0.062 & 0.063 & 0.064 & 0.065 & 0.066 & 0.063 \\ 
& & $\widehat{\epsilon}_{\text{ELW2S}}$ & 0.052 & 0.050 & 0.053 & 0.059 & 0.065 & 0.063 & 0.056 & 0.057 & 0.057 \\ 
\cmidrule{2-12}
& Variance & $\widehat{\epsilon}_{\text{aggvar}}$ & 0.017 & 0.018 & 0.018 & 0.018 & 0.019 & 0.018 & 0.018 & 0.017 & 0.018 \\ 
& &  $\widehat{\epsilon}_{\text{diffvar}}$ & 0.022 & 0.022 & 0.022 & 0.022 & 0.022 & 0.023 & 0.024 & 0.024 & 0.023 \\ 
& & $\widehat{\epsilon}_{\text{absval}}$ & 0.017 & 0.018 & 0.018 & 0.019 & 0.019 & 0.019 & 0.019 & 0.017 & 0.018 \\ 
& & $\widehat{\epsilon}_{\text{Higuchi}}$ & 0.010 & 0.011 & 0.012 & 0.013 & 0.014 & 0.014 & 0.014 & 0.014 & 0.013 \\ 
& & $\widehat{\epsilon}_{\text{Peng}}$ & 0.003 & 0.003 & 0.003 & 0.004 & 0.004 & 0.004 & 0.004 & 0.005 & 0.004 \\ 
 & & $\widehat{\epsilon}_{\text{RS}}$ & \textBF{0.001} & \textBF{0.002} & \textBF{0.002} & \textBF{0.002} & \textBF{0.002} & \textBF{0.002} & \textBF{0.002} & \textBF{0.002} & \textBF{0.002} \\ 
& & $\widehat{\epsilon}_{\text{RAR}}$ & 0.008 & 0.008 & 0.008 & 0.008 & 0.008 & 0.008 & 0.008 & 0.008 & 0.008 \\ 
& & $\widehat{\epsilon}_{\text{per}}$ & 0.014 & 0.014 & 0.014 & 0.014 & 0.014 & 0.014 & 0.014 & 0.014 & 0.014 \\ 
& & $\widehat{\epsilon}_{\text{boxper}}$ & 0.003 & 0.003 & 0.003 & 0.003 & 0.003 & 0.003 & 0.004 & 0.004 & 0.003 \\ 
& & $\widehat{\epsilon}_{\text{GPH}}$ & 0.029 & 0.029 & 0.028 & 0.028 & 0.028 & 0.028 & 0.029 & 0.028 & 0.028 \\ 
& & $\widehat{\epsilon}_{\text{SGPH}}$ & 0.016 & 0.017 & 0.017 & 0.017 & 0.017 & 0.018 & 0.018 & 0.018 & 0.017 \\ 
& & $\widehat{\epsilon}_{\text{Wavelet}}$ & 0.027 & 0.027 & 0.027 & 0.027 & 0.027 & 0.026 & 0.026 & 0.026 & 0.027 \\ 
& & $\widehat{\epsilon}_{\text{Local\_W}}$ & 0.006 & 0.006 & 0.006 & 0.006 & 0.006 & 0.006 & 0.006 & 0.005 & 0.006 \\ 
& & $\widehat{\epsilon}_{\text{Local\_W\_T}}$ & 0.023 & 0.023 & 0.022 & 0.022 & 0.020 & 0.019 & 0.017 & 0.013 & 0.020 \\ 
& & $\widehat{\epsilon}_{\text{Hou\_Perron}}$ & 0.008 & 0.008 & 0.008 & 0.008 & 0.009 & 0.011 & 0.012 & 0.013 & 0.010 \\ 
& & $\widehat{\epsilon}_{\text{ELW}}$ & 0.004 & 0.004 & 0.004 & 0.004 & 0.004 & 0.004 & 0.004 & 0.004 & 0.004 \\ 
& & $\widehat{\epsilon}_{\text{ELW2S}}$ & 0.006 & 0.007 & 0.009 & 0.009 & 0.008 & 0.006 & 0.006 & 0.007 & 0.007 \\   \cmidrule{2-12}
& MSE & $\widehat{\epsilon}_{\text{aggvar}}$ & 0.021 & 0.024 & 0.027 & 0.030 & 0.033 & 0.037 & 0.041 & 0.046 & 0.032 \\ 
& &  $\widehat{\epsilon}_{\text{diffvar}}$ & 0.074 & 0.071 & 0.072 & 0.072 & 0.068 & 0.069 & 0.069 & 0.070 & 0.071 \\ 
& & $\widehat{\epsilon}_{\text{absval}}$ & 0.017 & 0.018 & 0.019 & 0.019 & 0.020 & 0.021 & 0.023 & 0.025 & 0.020 \\ 
& & $\widehat{\epsilon}_{\text{Higuchi}}$ & 0.011 & 0.013 & 0.014 & 0.016 & 0.018 & 0.021 & 0.023 & 0.027 & 0.018 \\ 
& & $\widehat{\epsilon}_{\text{Peng}}$ & 0.010 & 0.008 & 0.007 & 0.007 & 0.006 & \textBF{0.006} & \textBF{0.005} & 0.005 & 0.007 \\ 
 & & $\widehat{\epsilon}_{\text{RS}}$ & \textBF{0.006} & \textBF{0.003} & \textBF{0.002} & \textBF{0.002} & \textBF{0.005} & 0.010 & 0.016 & 0.026 & 0.009 \\ 
& & $\widehat{\epsilon}_{\text{RAR}}$ & 0.023 & 0.017 & 0.013 & 0.010 & 0.008 & 0.008 & 0.010 & 0.013 & 0.013 \\ 
& & $\widehat{\epsilon}_{\text{per}}$ & 0.015 & 0.016 & 0.016 & 0.017 & 0.018 & 0.019 & 0.020 & 0.022 & 0.018 \\ 
& & $\widehat{\epsilon}_{\text{boxper}}$ & 0.085 & 0.085 & 0.085 & 0.085 & 0.085 & 0.085 & 0.085 & 0.085 & 0.085 \\ 
& & $\widehat{\epsilon}_{\text{GPH}}$ & 0.029 & 0.029 & 0.028 & 0.029 & 0.029 & 0.029 & 0.029 & 0.029 & 0.029 \\ 
& & $\widehat{\epsilon}_{\text{SGPH}}$ & 0.017 & 0.018 & 0.018 & 0.018 & 0.018 & 0.018 & 0.018 & 0.018 & 0.018 \\ 
& & $\widehat{\epsilon}_{\text{Wavelet}}$ & 0.038 & 0.037 & 0.036 & 0.035 & 0.034 & 0.033 & 0.032 & 0.032 & 0.035 \\ 
& & $\widehat{\epsilon}_{\text{Local\_W}}$ & 0.006 & 0.006 & 0.006 & 0.006 & 0.006 & 0.006 & 0.006 & \textBF{0.005} & \textBF{0.006} \\ 
& & $\widehat{\epsilon}_{\text{Local\_W\_T}}$ & 0.023 & 0.023 & 0.023 & 0.022 & 0.021 & 0.019 & 0.017 & 0.013 & 0.020 \\ 
& & $\widehat{\epsilon}_{\text{Hou\_Perron}}$ & 0.008 & 0.008 & 0.008 & 0.008 & 0.009 & 0.011 & 0.012 & 0.013 & 0.010 \\ 
& & $\widehat{\epsilon}_{\text{ELW}}$ & 0.008 & 0.008 & 0.008 & 0.008 & 0.008 & 0.008 & 0.008 & 0.008 & 0.008 \\ 
& & $\widehat{\epsilon}_{\text{ELW2S}}$ & 0.009 & 0.010 & 0.011 & 0.012 & 0.012 & 0.010 & 0.009 & 0.010 & 0.010 \\ 
\cmidrule{2-12}
1000 & Bias & $\widehat{\epsilon}_{\text{aggvar}}$ & -0.047 & -0.057 & -0.067 & -0.078 & -0.091 & -0.106 & -0.123 & -0.141 & -0.089 \\ 
& &  $\widehat{\epsilon}_{\text{diffvar}}$ & 0.201 & 0.199 & 0.200 & 0.199 & 0.194 & 0.194 & 0.193 & 0.194 & 0.197 \\ 
& & $\widehat{\epsilon}_{\text{absval}}$ & 0.016 & 0.006 & -0.004 & -0.015 & -0.028 & -0.043 & -0.060 & -0.079 & -0.026 \\ 
& & $\widehat{\epsilon}_{\text{Higuchi}}$ & -0.011 & -0.017 & -0.024 & -0.031 & -0.040 & -0.051 & -0.065 & -0.083 & -0.040 \\ 
& & $\widehat{\epsilon}_{\text{Peng}}$ & 0.065 & 0.057 & 0.049 & 0.042 & 0.035 & 0.029 & 0.024 & 0.019 & 0.040 \\ 
 & & $\widehat{\epsilon}_{\text{RS}}$ & 0.060 & 0.033 & 0.006 & -0.022 & -0.051 & -0.080 & -0.110 & -0.142 & -0.038 \\ 
& & $\widehat{\epsilon}_{\text{RAR}}$ & 0.112 & 0.089 & 0.065 & 0.040 & 0.014 & -0.012 & -0.040 & -0.068 & 0.025 \\ 
& & $\widehat{\epsilon}_{\text{per}}$ & 0.029 & 0.033 & 0.038 & 0.042 & 0.047 & 0.052 & 0.058 & 0.064 & 0.045 \\ 
& & $\widehat{\epsilon}_{\text{boxper}}$ & 0.191 & 0.191 & 0.192 & 0.192 & 0.193 & 0.193 & 0.194 & 0.195 & 0.193 \\ 
& & $\widehat{\epsilon}_{\text{GPH}}$ & 0.006 & 0.007 & 0.008 & 0.011 & 0.013 & 0.015 & 0.019 & 0.023 & 0.013 \\ 
& & $\widehat{\epsilon}_{\text{SGPH}}$ & -0.023 & -0.024 & -0.024 & -0.023 & -0.022 & -0.020 & -0.017 & -0.013 & -0.021 \\ 
& & $\widehat{\epsilon}_{\text{Wavelet}}$ & 0.077 & 0.073 & 0.069 & 0.066 & 0.062 & 0.059 & 0.056 & 0.053 & 0.064 \\ 
& & $\widehat{\epsilon}_{\text{Local\_W}}$ & -0.005 & -0.004 & -0.004 & -0.003 & -0.002 & -0.001 & \textBF{0.001} & 0.002 & -0.002 \\ 
& & $\widehat{\epsilon}_{\text{Local\_W\_T}}$ & 0.010 & 0.010 & 0.010 & 0.010 & 0.010 & 0.008 & 0.005 & \textBF{-0.002} & 0.008 \\ 
& & $\widehat{\epsilon}_{\text{Hou\_Perron}}$ & \textBF{0.003} & \textBF{0.003} & \textBF{0.002} & \textBF{0.001} & \textBF{0.001} & \textBF{-0.000} & -0.001 & -0.002 & \textBF{0.001} \\ 
& & $\widehat{\epsilon}_{\text{ELW}}$ & 0.040 & 0.040 & 0.041 & 0.041 & 0.042 & 0.042 & 0.043 & 0.044 & 0.042 \\ 
& & $\widehat{\epsilon}_{\text{ELW2S}}$ & 0.034 & 0.029 & 0.029 & 0.041 & 0.047 & 0.046 & 0.040 & 0.037 & 0.038 \\  \cmidrule{2-12}
& Variance & $\widehat{\epsilon}_{\text{aggvar}}$ & 0.011 & 0.012 & 0.012 & 0.012 & 0.012 & 0.012 & 0.012 & 0.011 & 0.012 \\ 
& &  $\widehat{\epsilon}_{\text{diffvar}}$ & 0.011 & 0.011 & 0.011 & 0.011 & 0.012 & 0.013 & 0.012 & 0.013 & 0.012 \\ 
& & $\widehat{\epsilon}_{\text{absval}}$ & 0.011 & 0.012 & 0.012 & 0.012 & 0.012 & 0.012 & 0.012 & 0.012 & 0.012 \\ 
& & $\widehat{\epsilon}_{\text{Higuchi}}$ & 0.005 & 0.005 & 0.006 & 0.006 & 0.007 & 0.007 & 0.007 & 0.007 & 0.006 \\ 
& & $\widehat{\epsilon}_{\text{Peng}}$ & 0.002 & 0.002 & 0.002 & 0.002 & 0.002 & 0.003 & 0.003 & 0.003 & 0.002 \\ 
 & & $\widehat{\epsilon}_{\text{RS}}$ & \textBF{0.001} & \textBF{0.001} & \textBF{0.001} & \textBF{0.001} & \textBF{0.001} & \textBF{0.002} & \textBF{0.002} & \textBF{0.002} & \textBF{0.001} \\ 
& & $\widehat{\epsilon}_{\text{RAR}}$ & 0.006 & 0.006 & 0.006 & 0.006 & 0.006 & 0.006 & 0.006 & 0.006 & 0.006 \\ 
& & $\widehat{\epsilon}_{\text{per}}$ & 0.006 & 0.006 & 0.006 & 0.006 & 0.006 & 0.006 & 0.006 & 0.006 & 0.006 \\ 
& & $\widehat{\epsilon}_{\text{boxper}}$ & 0.002 & 0.002 & 0.002 & 0.002 & 0.002 & 0.002 & 0.002 & 0.002 & 0.002 \\ 
& & $\widehat{\epsilon}_{\text{GPH}}$ & 0.019 & 0.019 & 0.020 & 0.019 & 0.019 & 0.020 & 0.020 & 0.020 & 0.020 \\ 
& & $\widehat{\epsilon}_{\text{SGPH}}$ & 0.011 & 0.011 & 0.012 & 0.012 & 0.012 & 0.012 & 0.012 & 0.013 & 0.012 \\ 
& & $\widehat{\epsilon}_{\text{Wavelet}}$ & 0.012 & 0.012 & 0.012 & 0.012 & 0.012 & 0.012 & 0.012 & 0.012 & 0.012 \\ 
& & $\widehat{\epsilon}_{\text{Local\_W}}$ & 0.003 & 0.003 & 0.003 & 0.003 & 0.003 & 0.003 & 0.003 & 0.003 & 0.003 \\ 
& & $\widehat{\epsilon}_{\text{Local\_W\_T}}$ & 0.014 & 0.014 & 0.014 & 0.014 & 0.013 & 0.013 & 0.011 & 0.009 & 0.013 \\ 
& & $\widehat{\epsilon}_{\text{Hou\_Perron}}$ & 0.004 & 0.004 & 0.004 & 0.005 & 0.005 & 0.005 & 0.005 & 0.006 & 0.005 \\ 
& & $\widehat{\epsilon}_{\text{ELW}}$ & 0.002 & 0.002 & 0.002 & 0.002 & 0.002 & 0.002 & 0.002 & 0.002 & 0.002 \\ 
& & $\widehat{\epsilon}_{\text{ELW2S}}$ & 0.003 & 0.004 & 0.005 & 0.006 & 0.005 & 0.004 & 0.003 & 0.004 & 0.004 \\\cmidrule{2-12}
& MSE & $\widehat{\epsilon}_{\text{aggvar}}$ & 0.013 & 0.015 & 0.016 & 0.018 & 0.020 & 0.023 & 0.027 & 0.031 & 0.020 \\ 
& &  $\widehat{\epsilon}_{\text{diffvar}}$ & 0.052 & 0.050 & 0.051 & 0.051 & 0.050 & 0.050 & 0.050 & 0.051 & 0.051 \\ 
& & $\widehat{\epsilon}_{\text{absval}}$ & 0.012 & 0.012 & 0.012 & 0.012 & 0.013 & 0.014 & 0.016 & 0.018 & 0.014 \\ 
& & $\widehat{\epsilon}_{\text{Higuchi}}$ & 0.005 & 0.006 & 0.006 & 0.007 & 0.008 & 0.010 & 0.011 & 0.014 & 0.008 \\ 
& & $\widehat{\epsilon}_{\text{Peng}}$ & 0.006 & 0.005 & 0.005 & 0.004 & 0.004 & 0.003 & \textBF{0.003} & 0.003 & 0.004 \\ 
 & & $\widehat{\epsilon}_{\text{RS}}$ & 0.005 & \textBF{0.002} & \textBF{0.001} & \textBF{0.002} & 0.004 & 0.008 & 0.014 & 0.022 & 0.007 \\ 
& & $\widehat{\epsilon}_{\text{RAR}}$ & 0.019 & 0.014 & 0.011 & 0.008 & 0.007 & 0.007 & 0.008 & 0.011 & 0.010 \\ 
& & $\widehat{\epsilon}_{\text{per}}$ & 0.007 & 0.007 & 0.007 & 0.008 & 0.008 & 0.009 & 0.009 & 0.010 & 0.008 \\ 
& & $\widehat{\epsilon}_{\text{boxper}}$ & 0.038 & 0.038 & 0.038 & 0.039 & 0.039 & 0.039 & 0.039 & 0.040 & 0.039 \\ 
& & $\widehat{\epsilon}_{\text{GPH}}$  & 0.019 & 0.019 & 0.020 & 0.019 & 0.019 & 0.020 & 0.020 & 0.020 & 0.020 \\ 
& & $\widehat{\epsilon}_{\text{SGPH}}$ & 0.012 & 0.012 & 0.012 & 0.012 & 0.012 & 0.012 & 0.013 & 0.013 & 0.012 \\ 
& & $\widehat{\epsilon}_{\text{Wavelet}}$ & 0.017 & 0.017 & 0.017 & 0.016 & 0.016 & 0.015 & 0.015 & 0.015 & 0.016 \\ 
& & $\widehat{\epsilon}_{\text{Local\_W}}$ & \textBF{0.003} & 0.003 & 0.003 & 0.003 & \textBF{0.003} & \textBF{0.003} & 0.003 & \textBF{0.003} & \textBF{0.003} \\ 
& & $\widehat{\epsilon}_{\text{Local\_W\_T}}$ & 0.014 & 0.014 & 0.014 & 0.014 & 0.013 & 0.013 & 0.011 & 0.009 & 0.013 \\ 
& & $\widehat{\epsilon}_{\text{Hou\_Perron}}$ & 0.004 & 0.004 & 0.004 & 0.005 & 0.005 & 0.005 & 0.005 & 0.006 & 0.005 \\ 
& & $\widehat{\epsilon}_{\text{ELW}}$ & 0.004 & 0.004 & 0.004 & 0.004 & 0.004 & 0.004 & 0.004 & 0.004 & 0.004 \\ 
& & $\widehat{\epsilon}_{\text{ELW2S}}$ & 0.004 & 0.004 & 0.006 & 0.008 & 0.008 & 0.006 & 0.005 & 0.005 & 0.006 \\
\bottomrule
\end{longtable}
\end{center}

In the supplement, we also consider the same functional models with the long-run covariance estimated from \cite{RS17}. Despite two ways of estimating the long-run covariance function, the summary statistics of the estimation bias, variance and MSE differ marginally and our recommendation for the best estimator, in terms of bias, variance, and MSE, remains the same.

The ranking of estimators might depend on the particular short-memory dependence of the time series, so we have also explored some variants of Case 1 and Case 2 for a selected subset of estimates in Table~\ref{tab:111}. 

By altering the coefficient in the kernel function, it changes the temporal dependence structure from weak to strong dependence in the ARFIMA$(1, d, 0)$ case. In the ARFIMA$(1, d, 0)$ model, the coefficient 0.34 produces a $L_2$ norm of 0.5, implying a moderate temporal dependence. To have the $L_2$ norm of $\phi_1$ to be 0.1, we change that coefficient from 0.34 to 0.068. To have the $L_2$ norm of $\phi_1$ to be 0.9, we change that coefficient from 0.34 to 0.612. From the MSE of the functional ARFIMA$(1, d, 0)$ model in Table~\ref{tab:1}, the Peng estimator produces the smallest estimation error, followed by the local Whittle estimator. From Table~\ref{tab:111}, we observe that the Peng estimator still performs better than the local Whittle estimator under various degrees of dependence.

Similarly, by altering the coefficients in the kernel functions, it alters the temporal dependence structure from weak to strong dependence in the ARFIMA$(1, d, 1)$ case. In the functional ARFIMA$(1, d, 1)$ model, the coefficients $(0.34, 1.5)$ produces a $L_2$ norm of 0.5 for both the AR and MA components implying a moderate temporal dependence. To have the $L_2$ norm of $\phi_1$ and $\psi_1$ to be 0.1, we change the coefficients from $(0.34, 1.5)$ to $(0.068, 0.059)$. To have the $L_2$ norm of $\phi_1$ and $\psi_1$ to be 0.9, we change the coefficients from $(0.34, 1.5)$ to $(0.612, 4.765)$. From the MSE of the functional ARFIMA$(1, d, 0)$ model, in Table~\ref{tab:11}, the local Whittle estimator produces the smallest estimation error, followed by the Peng estimator. From Table~\ref{tab:111}, we observe that the local Whittle estimator still performs better than the Peng estimator under the strong dependence, but not so under the weak dependence. From this example, various degrees of short-memory temporal dependence can affect the estimation accuracy of the long-memory parameter.

In the supplement, we also consider the same functional models with the long-run covariance estimated from \cite{RS17}. Despite two ways of estimating the long-run covariance function, the summary statistics of the MSE differ marginally, and our recommendation for the best estimator, in terms of MSE, remains the same.

\begin{center}
\tabcolsep 0.125in
\begin{longtable}{@{}llllrrrrrrrrrr@{}}
\caption{With sample size $n=250, 500, 1000$, we evaluate and compare the finite-sample MSE ($\times 100$) of the estimation error $\Delta d$ under the functional ARFIMA($1,d,0$) and ARFIMA($1,d,1$) models, where the long-run covariance was estimated from \cite{LRS18}.}\label{tab:111}
\\
\toprule
$\rho$ & Estimator & $n$ &  \multicolumn{8}{c}{$d$} & Overall\\
 & & & 0.05 & 0.10 & 0.15 & 0.20 & 0.25 & 0.30 & 0.35 & 0.40 & \\
\endfirsthead
\toprule
$\rho$ & Estimator & $n$ & \multicolumn{8}{c}{$d$} & Overall \\
 & & & 0.05 & 0.10 & 0.15 & 0.20 & 0.25 & 0.30 & 0.35 & 0.40 & \\
\midrule
\endhead
\midrule
\multicolumn{12}{r}{Continued on next page} 
\endfoot
\endlastfoot
\midrule
\multicolumn{4}{l}{\hspace{-.13in}{\underline{DGP = ARFIMA($1, d, 0$)}}} \\
0.9 & Peng & 250 & 0.96 & 0.87 & 0.81 & 0.78 & 0.76 & 0.77 & 0.78 & 0.81 & 0.82 \\
	& 	& 500 & 0.53 & 0.47 & 0.43 & 0.41 & 0.41 & 0.42 & 0.43 & 0.45 & 0.44 \\
	&	& 1000 & 0.33 & 0.29 & 0.27 & 0.26 & 0.26 & 0.26 & 0.27 & 0.29 & 0.28 \\
\\
& Local\_W & 250 & 0.95 & 0.94 & 0.94 & 0.93 & 0.92 & 0.90 & 0.86 & 0.76 & 0.90 \\
	&	& 500 & 0.61 & 0.61 & 0.61 & 0.61 & 0.60 & 0.60 & 0.58 & 0.54 & 0.60 \\
	&	& 1000 & 0.35 & 0.35 & 0.35 & 0.35 & 0.35 & 0.34 & 0.34 & 0.33 & 0.35 \\
\\
0.1 & Peng & 250 & 0.50 & 0.56 & 0.63 & 0.70 & 0.76 & 0.83 & 0.89 & 0.94 & 0.73 \\
	&	& 500 & 0.27 & 0.31 & 0.35 & 0.39 & 0.44 & 0.48 & 0.52 & 0.55 & 0.41 \\
	&	& 1000 & 0.17 & 0.19 & 0.22 & 0.25 & 0.28 & 0.30 & 0.33 & 0.35 & 0.26 \\
\\
& Local\_W & 250 & 0.92 & 0.91 & 0.89 & 0.86 & 0.83 & 0.80 & 0.76 & 0.68 & 0.83  \\
	&	& 500 & 0.63 & 0.62 & 0.62 & 0.61 & 0.60 & 0.59 & 0.56 & 0.51 & 0.59 \\
	&	& 1000 & 0.36 & 0.36 & 0.35 & 0.35 & 0.35 & 0.34 & 0.34 & 0.32 & 0.35 \\	
\\
\multicolumn{4}{l}{\hspace{-.13in}{\underline{DGP = ARFIMA($1, d, 1$)}}} \\
0.9 &  Local\_W & 250 & 0.96 & 0.96 & 0.96 & 0.96 & 0.96 & 0.95 & 0.89 & 0.73 & 0.92 \\
	&	& 500 &  0.59 & 0.58 & 0.58 & 0.58 & 0.59 & 0.58 & 0.57 & 0.51 & 0.57 \\
	&	& 1000 & 0.34 & 0.34 & 0.34 & 0.34 & 0.34 & 0.34 & 0.34 & 0.32 & 0.34 \\
\\
& Peng & 250 & 1.80 & 1.58 & 1.40 & 1.26 & 1.15 & 1.08 & 1.02 & 0.99 & 1.29 \\
	&	& 500 & 1.05 & 0.90 & 0.79 & 0.70 & 0.64 & 0.59 & 0.56 & 0.55 & 0.72	\\
	&	& 1000 & 0.65 & 0.55 & 0.48 & 0.42 & 0.38 & 0.36 & 0.34 & 0.33 & 0.44	\\
\\	
0.1 & Local\_W	& 250 & 0.92 & 0.90 & 0.89 & 0.86 & 0.85 & 0.80 & 0.76 & 0.68 & 0.83 \\
	&	& 500 & 0.62 & 0.62 & 0.61 & 0.61 & 0.60 & 0.59 & 0.56 & 0.51 & 0.59 \\
	&	& 1000 & 0.36 & 0.35 & 0.35 & 0.35 & 0.34 & 0.34 & 0.34 & 0.32 & 0.34 \\
\\
& Peng & 250 & 0.50 & 0.55 & 0.61 & 0.67 & 0.73 & 0.80 & 0.85 & 0.90 & 0.70 \\
	&	& 500 & 0.27 &	0.30 & 0.33 & 0.38 & 0.42 & 0.46 & 0.50 & 0.53 & 0.40 \\
	&	& 1000 & 0.17 & 0.19 & 0.21 & 0.24 & 0.26 & 0.29 & 0.31 & 0.33 & 0.25	\\	
\bottomrule
\end{longtable}
\end{center}

\section{Conclusion}

Through a series of simulation studies, we identify the most accurate estimator for estimating the long-memory parameter in a functional ARFIMA model. For a functional ARFIMA$(1, d, 0)$ with various $d$ values, the local Whittle estimator with tapering produces the smallest bias; the rescaled range estimator produces the smallest variance, and the Peng estimator produces the smallest mean square error. For a functional ARFIMA$(1, d, 1)$ with various $d$ values, the local Whittle estimator produces the smallest bias with $n=250$ and the Hou-Perron estimator produces the smallest bias with $n=500$ and $n=1000$; the rescaled range estimator produces the smallest variance; and the local Whittle estimator produces the smallest mean square error.

For the functional ARFIMA$(1, d, 0)$ model, the Peng estimator consistently produces the smallest estimation error under weak, moderate, and strong dependence. For the functional ARFIMA$(1, d, 1)$ model, the local Whittle estimator produces the smallest estimation errors for moderate and strong dependence, while the Peng estimator produces the smallest estimation error for weak dependence. In summary, various degrees of short-memory temporal dependence can affect the estimation accuracy of the long-memory parameter.

There are several ways in which the present study can be further extended, and we briefly outline three:
\begin{inparaenum}
\item[1)] Estimation of a time-varying long-memory parameter. 
\item[2)] Consider a wavelet-based or Fourier-based multivariate Whittle estimation for a multivariate time series of principal component scores. The multivariate estimation method should be more efficient than a univariate estimation method, subject to each set of principal component scores shows a similar degree of persistence. Should this condition fails to satisfy, it may lead to biased estimates. 
\item[3)] Propose tests for detecting the presence of long-memory.
\end{inparaenum}

%\section*{Acknowledgment}

%The author thanks the financial support of a faculty grant and computing facility provided at the College of Business and Economics at the Australian National University. 

%\newpage
\bibliographystyle{agsm}
\bibliography{long_memory.bib}

\end{document}